\documentclass[11pt, leqno]{article} 
\usepackage{latexsym, epsfig, amssymb, amsmath}
\usepackage[dvips]{color}
\usepackage{graphicx}
\usepackage{anysize}

\makeatletter

\renewcommand{\baselinestretch}{1.3} 
\marginsize{1.25 in}{1.25 in}{0.8 in}{1.4 in}

\def\singlespace{\def\baselinestretch{1.3}\@normalsize}

\newcounter{sect}
\newcommand{\sect}{\refstepcounter{sect} \noindent\quad{\bf \arabic{sect}.\ }}
\newcounter{subsect}
\newcommand{\subsect}{\refstepcounter{subsect} \noindent\quad{\arabic{sect}.\arabic{subsect}.\ }}
\newcounter{subsubsect}

\newcounter{lem}

\newcounter{prop}

\newcounter{thm}

\renewcommand{\theequation}{\arabic{sect}.\arabic{equation}}

\makeatother

\newcommand{\bb}{\mbox{\bf b}}
\newcommand{\bff}{\mbox{\bf f}}

\newcommand{\by}{\mbox{\bf y}}
\newcommand{\bA}{\mbox{\bf A}}
\newcommand{\bB}{\mbox{\bf B}}
\newcommand{\bC}{\mbox{\bf C}}
\newcommand{\bD}{\mbox{\bf D}}
\newcommand{\bE}{\mbox{\bf E}}
\newcommand{\bF}{\mbox{\bf F}}
\newcommand{\bG}{\mbox{\bf G}}
\newcommand{\bH}{\mbox{\bf H}}
\newcommand{\bQ}{\mbox{\bf Q}}
\newcommand{\bX}{\mbox{\bf X}}
\newcommand{\bY}{\mbox{\bf Y}}
\newcommand{\bone}{\mbox{\bf 1}}
\newcommand{\bzero}{\mbox{\bf 0}}
\newcommand{\bveps}{\mbox{\boldmath $\varepsilon$}}
\newcommand{\bet}{\mbox{\boldmath $\eta$}}
\newcommand{\bxi}{\mbox{\boldmath $\xi$}}
\newcommand{\bmu}{\mbox{\boldmath $\mu$}}

\newcommand{\hB}{\widehat \bB}
\newcommand{\hb}{\widehat \bb}
\newcommand{\hE}{\widehat \bE}
\newcommand{\hvar}{\widehat \var}
\newcommand{\hcov}{\widehat \cov}

\newcommand{\hSig}{\widehat\Sig}
\newcommand{\hsig}{\widehat\sigma}
\newcommand{\hmu}{\widehat\bmu}
\newcommand{\hxi}{\widehat\bxi}
\newcommand{\heq}{\ \widehat=\ }
\newcommand{\sam}{_{\text{sam}}}

\newcommand{\var}{\mathrm{var}}
\newcommand{\cov}{\mathrm{cov}}
\newcommand{\Sig}{\mathbf{\Sigma}}
\newcommand{\veps}{\varepsilon}
\newcommand{\tr}{\mathrm{tr}}
\newcommand{\diag}{\mathrm{diag}}
\newcommand{\vecc}{\mathrm{vec}}
\newcommand{\vech}{\mathrm{vech}}

\def\t {'}
\def\toD{\overset{\mathrm{D}}{\longrightarrow}}
\def\toP{\overset{\mathrm{P}}{\longrightarrow}}

\begin{document}

\title{\vspace{0in}\bf High Dimensional Covariance Matrix Estimation Using a
Factor Model
\thanks{Financial support from the NSF under grant DMS-0532370 is gratefully acknowledged.
Address for correspondence: Jinchi Lv, Department of Mathematics,
Princeton University, Princeton, NJ 08544. Phone: (609) 258-9433.
E-mail: jlv@princeton.edu.}
\date{August 12, 2006}
\smallskip
\author{\textsc{By Jianqing Fan, Yingying Fan and
Jinchi Lv}\\
\em Princeton University} }

\maketitle

\vspace{-0.15 in}
\begin{singlespace}
\begin{quotation}
High dimensionality comparable to sample size is common in many
statistical problems.
We examine covariance matrix estimation in the asymptotic framework
that the dimensionality $p$ tends to $\infty$ as the sample size $n$
increases. Motivated by the Arbitrage Pricing Theory in finance, a
multi-factor model is employed to reduce dimensionality and to
estimate the covariance matrix.
The factors are 
observable and the number of
factors $K$ is allowed to grow with 
$p$.  We
investigate impact of 
$p$ and $K$ on the performance of the model-based covariance matrix
estimator. Under mild assumptions, we have established 
convergence rates and asymptotic normality of the model-based
estimator.
Its performance is compared with that of the sample covariance
matrix.  
We identify situations under which the factor approach 
increases performance substantially or marginally. 
The impacts of covariance matrix estimation on portfolio allocation
and risk management 
are studied. The asymptotic
results are supported by a thorough simulation study.
\end{quotation}
\end{singlespace}

\bigskip

 {\em Short Title}: Large Covariance Matrix
Estimation.



\bigskip
\newpage

\sect {\bf Introduction.}

\bigskip

\subsect {\em Background}.\quad Covariance matrix estimation is
fundamental for almost all areas of multivariate analysis and many
other applied problems. In particular, covariance matrices and their
inverses play a central role in risk management and portfolio
allocation. For example, the smallest and largest eigenvalues of a
covariance matrix are related to the minimum and maximum variances
of the selected portfolio, respectively, and the eigenvectors are
related to portfolio allocation. Therefore, we need a good
covariance matrix estimator inverting which does not excessively
amplify the estimation error. See Goldfarb and Iyengar (2003) for
applications of covariance matrices to portfolio selections and
Johnstone (2001) for their statistical implications.

Estimating high-dimensional covariance matrices is intrinsically
challenging. For example, in portfolio allocation and risk
management, the number of stocks $p$, which is typically of the same
order as the sample size $n$, can well be in the order of hundreds.
In particular, when $p=200$ there are more than 20,000 parameters in
the covariance matrix. Yet, the available sample size is usually in
the order of hundreds or a few thousands because longer time series
(larger $n$) increases modeling bias. For instance, by taking daily
data of the past three years we have only roughly $n=750$. So it is
hard or even unrealistic to estimate covariance matrices without
imposing any structure (see the rejoinder in Fan, 2005).

Factor models have been widely used both theoretically and
empirically in economics and finance. Derived by Ross (1976, 1977)
using the Arbitrage Pricing Theory (APT) and by Chamberlain and
Rothschild (1983) in a large economy, the multi-factor model states
that the excessive return of any asset $Y_i$ over the risk-free
interest rate satisfies
\begin{equation} \label{100}
Y_i=b_{i1}f_1+\cdots+b_{iK}f_K+\varepsilon_i,\quad i=1,\cdots,p,
\end{equation}
where $f_1,\cdots,f_K$ are the excessive returns of $K$ factors,
$b_{ij}$, $i=1,\cdots,p$, $j=1,\cdots,K$, are unknown factor
loadings, and $\veps_1,\cdots,\veps_p$ are $p$ idiosyncratic errors
uncorrelated given $f_1,\cdots,f_K$. In economics and finance
literature, factors are implicitly assumed to be observable and
there is a large literature contributed to construction of factors
(e.g. Fama and French, 1992, 1993).  The factor models have been
widely applied in economics and finance. See, for example, Ross
(1976, 1977), Engle and Watson (1981), Chamberlain (1983),
Chamberlain and Rothschild (1983), Diebold and Nerlove (1989), Fama
and French (1992, 1993), Aguilar and West (2000), and Stock and
Watson (2005) and references therein.  These are extensions of the
famous Capital Asset Pricing Model (CAPM) and can be regarded as
efforts to approximate the market portfolio in the CAPM.

Thanks to the multi-factor model (\ref{100}), if a few factors can
completely capture the cross-sectional risks, the number of
parameters in covariance matrix estimation can be significantly
reduced. For example, using the Fama-French three-factor model [Fama
and French (1992, 1993)], there are $4p$ instead of $p(p+1)/2$
parameters to be estimated. Despite the popularity of factor models
in the literature, the impact of dimensionality on the estimation
errors of covariance matrices and its applications to portfolio
allocation and risk management are poorly understood, so in this
paper, determined efforts are made on such an investigation. To make
the multi-factor model more realistic, we allow $K$ to grow with the
number of assets $p$ and hence with the sample size $n$. As a
result, we also investigate the impact of the number of factors on
the estimation of covariance matrices, as well as its applications
to portfolio allocation and risk management.  To appreciate the
derived rates of convergence, we compare them with those without
using the factor structure.  One natural candidate is the sample
covariance matrix. This also allows us to examine the impact of
dimensionality on the performance of the sample covariance matrix.
Our results can also be regarded as an important step to understand
the performance of factor models with unobservable factors.

The factor model has been extensively studied in the literature
[see, e.g. Scott (1966) and (1969), Browne (1987), Browne and
Shapiro (1987), and Yuan and Bentler (1997)], but traditional work
assumes the sample size $n$ tends to infinity while the
dimensionality $p$ and the number of factors $K$ are fixed. There is
a relatively small literature on studies of models with a diverging
number of parameters. See, for example, Huber (1973), Yohai and
Maronna (1979), Portnoy (1984, 1985), and Bai (2003). In particular,
Fan and Peng (2004) establish some asymptotic properties, as well as
an oracle property, for nonconcave penalized likelihood estimators
in the presence of a diverging number of parameters. One can further
refer to seminal reviews by Donoho (2000) and Fan and Li (2006) for
challenges of high dimensionality. But it still remains open to
examine factor models with diverging dimensionality and growing
number of factors for the purpose of covariance matrix estimation.

The traditional covariance matrix estimator, the sample covariance
matrix, is known to be unbiased, and it is invertible when the
dimensionality is no larger than the sample size. See, for example,
Eaton and Tyler (1991, 1994) for the asymptotic spectral
distributions of random matrices including sample covariance
matrices and their statistical implications. In the absence of prior
information about the population covariance matrix, the sample
covariance matrix is certainly a natural candidate in the case of
small dimensionality, but no longer performs very well for moderate
or large dimensionality [see, e.g. Lin and Perlman (1985) and
Johnstone (2001)]. Many approaches were proposed in the literature
to construct good covariance matrix estimators. Among them, two main
directions were taken. One is to remedy the sample covariance matrix
and construct a better one by using approaches such as shrinkage and
the eigen-method, etc. See, for example, Ledoit and Wolf (2004) and
Stein (1975). The other one is to reduce dimensionality by imposing
some structure on the data. Many structures, such as sparsity,
compound symmetry, and the autoregressive model, are widely used.
Various approaches were taken to seek a balance between the bias and
variance of covariance matrix estimators. See, for example, Dempster
(1972), Leonard and Hsu (1992), Chiu, Leonard and Tsui (1996),
Diggle and Verbyla (1998), Pourahmadi (2000), Boik (2002), Smith and
Kohn (2002), Wong, Carter and Kohn (2003), Wu and Pourahmadi (2003),
Huang, Liu and Pourahmadi (2004), and Li and Gui (2005).

\bigskip

\subsect {\em Covariance matrix estimation}.\quad We always denote
by $n$ the sample size, by $p$ the dimensionality, and by
$f_1,\cdots,f_K$ the $K$ observable factors, where $p$ grows with
sample size $n$ and $K$ increases with dimensionality $p$. For ease
of presentation, we rewrite factor model (\ref{100}) in matrix form
\begin{equation} \label{101}
\by=\bB_n\bff+\bveps,
\end{equation}
where $\by=(Y_1,\cdots,Y_p)\t$, $\bB_n=(\bb_1,\cdots,\bb_p)\t$ with
$\bb_i=(b_{n,i1},\cdots,b_{n,iK})\t$, $i=1,\cdots,p$,
$\bff=(f_1,\cdots,f_K)\t $, and $\bveps=(\veps_1,\cdots,\veps_p)\t
$. Throughout we assume that $E(\bveps|\bff)=\bzero$ and
$\cov(\bveps|\bff)=\Sig_{n,0}$ is diagonal. For brevity of notation,
we suppress the first subscript $n$ in some situations where the
dependence on $n$ is self-evident.

Let $(\bff_1,\by_1),\cdots,(\bff_n,\by_n)$ be $n$ independent and
identically distributed (i.i.d.) samples of $(\bff,\by)$. We
introduce here some notation used throughout the paper. Let
\[ \Sig_n=\cov(\by),\ \bX=(\bff_1,\cdots,\bff_n),\
\bY=(\by_1,\cdots,\by_n)\text{ and }\bE=(\bveps_1,\cdots,\bveps_n).
\] Under model (\ref{101}), we have
\begin{equation} \label{102}
\Sig_n=\cov(\bB_n\bff)+\cov(\bveps)=\bB_n\cov(\bff)\bB_n\t
+\Sig_{n,0}.
\end{equation}
A natural idea for estimating $\Sig_n$ is to plug in the
least-squares estimators of $\bB_n$, $\cov(\bff)$, and $\Sig_{n,0}$.
Therefore, we have a substitution estimator
\begin{equation} \label{105}
\hSig_n=\hB_n\hcov(\bff)\hB_n\t +\hSig_{n,0},
\end{equation}
where $\hB_n=\bY\bX\t (\bX\bX\t )^{-1}$ is the matrix of estimated
regression coefficients, $\hcov(\bff)=(n-1)^{-1}\bX\bX\t
-\{n(n-1)\}^{-1}\bX\bone\bone\t \bX\t$ is the sample covariance
matrix of the factors $\bff$, and
\[ \hSig_{n,0}=\diag\left(n^{-1}\hE\hE\t\right) \] is the diagonal matrix of
$n^{-1}\hE\hE\t$ with $\hE=\bY-\hB\bX$ the matrix of residuals. If
the factor model is not employed, then we have the sample covariance
matrix estimator
\begin{equation} \label{106}
\hSig\sam=\left(n-1\right)^{-1}\bY\bY\t
-\left\{n\left(n-1\right)\right\}^{-1}\bY\bone\bone\t \bY\t.
\end{equation}

This paper mainly provides a theoretical understanding of the factor
model with a diverging dimensionality and growing number of factors
for the purpose of covariance matrix estimation; it does not aim to
compare with other popular estimators. Throughout the paper, we
always contrast the performance of the covariance matrix estimator
$\hSig$ in (\ref{105}) with that of the sample covariance matrix
$\hSig\sam$ in (\ref{106}). With prior information of the true
factor structure, the substitution estimator $\hSig$ is expected to
perform better than $\hSig\sam$. However, this has not formally been
shown, especially when $p\rightarrow\infty$ and
$K\rightarrow\infty$, and this is not always true. In addition,
exact properties of this kind are not well understood. As the
problem is important for portfolio management, determined efforts
are devoted in regard to this.   Our conclusion can be summarized as
follows.
\begin{itemize}
\item $\hSig$ is always invertible, even if $p>n$, while
$\hSig\sam$ suffers from the problem of possibly being singular when
dimensionality $p$ is close to or larger than sample size $n$.

\item  The advantage of the factor model lies in the estimation of
the inverse of the covariance matrix, not the estimation of the
covariance matrix itself. When the parameters involve the inverse of
the covariance matrix, the factor model shows substantial gains,
whereas when the parameters involved the covariance matrix directly,
the factor model does not have much advantage. The latter is a
surprise to the conventional wisdom.

\item  Portfolio allocations involve the inverse of the covariance
matrix and the factor-model based estimates gain substantially,
whereas the risk management involves directly the covariance matrix
and the gain is only marginally.

\item $\hSig$ has asymptotic normality, while in general $\hSig\sam$
may not have asymptotic normality of the same kind.
\end{itemize}
These properties will be demonstrated in our paper as follows.

\bigskip

\subsect {\em Outline of the paper}.\quad In section 2 we discuss
some basic assumptions and present the sampling properties of the
estimator $\hSig$, as well as those of $\hSig\sam$. We study the
impacts of the covariance matrix estimation on portfolio allocation
and risk management in Section 3. A simulation study is presented in
Section 4, which augments our theoretical study. Section 5 contains
some concluding remarks. The proofs of our results are given in
Section 6. All the technical lemmas are relegated to the Appendix.

\bigskip
\setcounter{equation}{0} \setcounter{subsect}{0}
\setcounter{subsubsect}{0}

\sect {\bf Sampling properties.}\quad In this section we study the
sampling properties of $\hSig$ and $\hSig\sam$ with growing
dimensionality and number of factors. We discuss some basic
assumptions in Section 2.1. The sampling properties are presented in
Section 2.2.

In the presence of diverging dimensionality, we should carefully
choose appropriate norms for high dimensional matrices in different
situations. We first introduce some notation. We always denote by
$\lambda_1(\bA),\cdots,\lambda_q(\bA)$ the $q$ eigenvalues of a
$q\times q$ symmetric matrix $\bA$ in decreasing order. For any
matrix $\bA=(a_{ij})$,  its Frobenius norm is given by
\begin{equation} \label{107}
\|\bA\|=\left\{\tr(\bA\bA\t)\right\}^{1/2}.
\end{equation}
In particular, if $\bA$ is a $q\times q$ symmetric matrix, then
$\|\bA\|=\left\{\sum_{i=1}^q\lambda_i(\bA)^2\right\}^{1/2}$. The
Frobenius norm as well as many other matrix norms [see Horn and
Johnson (1985)] is intrinsically related to the eigenvalues or
singular values of matrices.

Despite its popularity, the Frobenius norm is not appropriate for
understanding the performance of the factor-model based estimation
of the covariance matrix. To see this, let us consider a simple
example. Suppose we know ideally that $\bB=\bone$ and
$\cov\left(\bveps|f\right)=I_{p}$ in model (\ref{101}) with a single
factor $f$. Then we have a substitution covariance matrix estimator
$\hSig=\bone\hvar(f)\bone\t+I_{p}$ as in (\ref{105}). It is a
classical result that
\[ E\left|\hvar(f)-\var(f)\right|^2=O(n^{-1}). \]
Thus by (\ref{102}), we have
\[ \hSig-\Sig=\bone\left[\hvar(f)-\var(f)\right]\bone\t \]
and the Frobenius norm
$\|\hSig-\Sig\|=\left|\hvar(f)-\var(f)\right|p$ picks up and
amplifies the estimation error from $\hvar(f)$. Consequently, \[
E\left\|\hSig-\Sig\right\|^2=O(n^{-1}p^2). \] On the other hand, by
assuming boundedness of the fourth moments of $\by$ across $n$, a
routine calculation reveals that
\[ E\left\|\hSig\sam-\Sig\right\|^2=O(n^{-1}p^2). \]
This shows that under Frobenius norm, $\hSig$ and $\hSig\sam$ have
the same convergence rate and perform roughly the same. Thus we
should seek other norms that fully employ the factor structure. By
assuming the eigenvalues of $\Sig$ are bounded away from 0 and
$\var(f)>0$, routine calculations show that
\[
\left\|\Sig^{-1/2}\left(\hSig\sam-\Sig\right)\Sig^{-1/2}\right\|=O_P(n^{-1/2}p^{3/2}),
\]
whereas $\|\Sig^{-1/2}(\hSig-\Sig)\Sig^{-1/2}\|=O_P(n^{-1/2})$.
Therefore, with prior information of the true factor structure,
$\hSig$ performs much better than $\hSig\sam$ from this point of
view.

Motivated by the above example, we first fix a sequence of positive
definite covariance matrices $\Sig_n$ of dimensionality $p_n$,
$n=1,2,\cdots$, and define a new norm
\begin{equation} \label{108}
\left\|\bA\right\|_{\Sig_n}=p_n^{-1/2}\left\|\Sig_n^{-1/2}\bA\Sig_n^{-1/2}\right\|
\end{equation}
for any $p_n\times p_n$ matrix $\bA$. In particular, we have
$\|\Sig_n\|_{\Sig_n}=p^{-1/2}\|I_{p}\|=1$. The inclusion of a
normalization factor $p^{-1/2}$ above is not essential and we
incorporate it to take into account the diverging dimensionality. As
seen below, under this new norm $\|\cdot\|_{\Sig}$, the consistency
rate in the factor approach is better than that in the sample
approach. Equivalently, we are investigating convergence rates under
the loss function
\begin{equation} \label{146}
L(\hSig,\Sig)=p^{1/2}\left\|\hSig-\Sig\right\|_{\Sig}=
\left\{\tr[\hSig\Sig^{-1}-I]^2\right\}^{1/2}.
\end{equation}
The above definition of the norm $\|\cdot\|_{\Sig}$ seems a bit
artificial and involves the inverse of the true covariance matrix,
but it is very similar to the entropy loss function proposed by
James and Stein (1961). See Section 4 for further details.
Intrinsically, this norm takes into account and fully employs the
factor structure. In fact, as shown in the above example, the
advantage of the factor structure lies in better performance of the
inverse $\hSig^{-1}$. We will see later in this section that
$\hSig^{-1}$ is a much better estimator of $\Sig^{-1}$ than
$\hSig\sam^{-1}$, and this advantage is carried further in portfolio
allocation.

\bigskip

\subsect {\em Some basic assumptions}.\quad Let $b_n=E\|\by\|^2$,
$c_n=\max_{1\leq i\leq K}E(f_i^4)$, and $d_n=\max_{1\leq i\leq
p}E(\veps_i^4)$.

\bigskip

\noindent(A) $(\bff_1,\by_1),\cdots,(\bff_n,\by_n)$ are $n$ i.i.d.
samples of $(\bff,\by)$. $E(\bveps|\bff)=\bzero$ and
$\cov(\bveps|\bff)=\Sig_{n,0}$ is diagonal. Also, the distribution
of $\bff$ is continuous and $K\leq p$.

\medskip

The first and second parts are usual conditions, and it is realistic
to put $K\leq p$. The assumption that $\bff$ has a continuous
distribution is made to ensure that the $K\times K$ matrix
$\bX\bX\t$ is invertible with probability one when $n\geq K$.
Clearly, the covariance matrix estimator $\hSig$ is positive
definite with probability one whenever $n\geq K$.  By the assumption
that the $K$ factors capture the cross-sectional risks, the
idiosyncratic noises are uncorrelated, so $\Sig_{n,0}$ is diagonal.

\medskip
\smallskip

\noindent(B) $b_n=O(p)$ and the sequences $c_n$ and $d_n$ are
bounded. Also, there exists a constant $\sigma_1>0$ such that
$\lambda_{K}(\cov(\bff))\geq\sigma_1$ for all $n$.

\medskip

This is a technical assumption. In view of
$E\|\by\|^2=\sum_{i=1}^pEy_i^2$, $b_n=O(p)$ is a reasonable
condition. The assumption $c_n=O(1)$ shows that the fourth moments
of $\bff$ are bounded across $n$, which facilitates the study of the
sample covariance matrix of $\bff$. The uniform lower bound imposed
on the eigenvalues of $\cov(\bff)$ helps the study of the inverse of
the sample covariance matrix of $\bff$ since $K\rightarrow\infty$,
and it along with $b_n=O(p)$ entails that $\|\bB_n\|=O(p^{1/2})$. It
is evident from our theoretical analysis that
$\lambda_{K}(\cov(\bff))$ can be allowed to tend to zero at some
rate, which results in slower convergence rates of the estimators.
But we do not pursue in this direction here.
\medskip
\smallskip

\noindent(C) There exists a constant $\sigma_2>0$ such that
$\lambda_{p}(\Sig_{n,0})\geq\sigma_2$ for all $n$.

\medskip

This is a reasonable assumption and ensures that all the eigenvalues
of $\Sig_n$'s are bounded away from 0 in view of (\ref{102}). In
particular, we have $\|\Sig_n^{-1}\|=O(p^{1/2})$.  Our theoretical
analysis applies to the case where $\lambda_{p}(\Sig_{n,0})$ tends
to zero at some rate, but we do not pursue along this direction for
simplicity.

\medskip
\smallskip

\noindent(D) The $K$ factors $f_1,\cdots,f_K$ are fixed across $n$,
and $p^{-1}\bB_n\t\bB_n\rightarrow \bA$ as $n\rightarrow\infty$ for
some $K\times K$ symmetric positive semidefinite matrix $\bA$.

\medskip

This assumption is used {\em only} to establish asymptotic normality
of the estimator $\hSig$, which facilitates statistical inferences.
In view of
$p^{-1}\bB_n\t\bB_n=p^{-1}(\bb_1\bb_1\t+\cdots+\bb_p\bb_p\t)$, this
assumption is reasonable when $K$ is fixed.

\bigskip

\subsect {\em Sampling properties}.\quad

\bigskip

\begin{thm} (Rates of convergence under Frobenius norm).\quad\em Under conditions (A) and
(B), we have $\left\|\hSig-\Sig\right\|=O_P(n^{-1/2}pK)$ and
$\left\|\hSig\sam-\Sig\right\|=O_P(n^{-1/2}pK)$. In addition, we
have
\[ \max_{1\leq k\leq
p}\left|\lambda_k(\hSig_n)-\lambda_k(\Sig_n)\right|=o_P\{(p^2K^2\log
n/n)^{1/2}\}
\] and
\[ \max_{1\leq k\leq
p}\left|\lambda_k(\hSig\sam)-\lambda_k(\Sig_n)\right|=o_P\{(p^2K^2\log
n/n)^{1/2}\}. \bigskip \]
\end{thm}
\smallskip \quad From this theorem, we see that under the Frobenius norm, the
dimensionality reduces rates of convergence by an order of $pK$,
which is the order of the number of parameters. The above rate of
eigenvalues of $\hSig$ is optimal. To see it, let us extend the
previous example by including $K$ factors $f_1,\cdots,f_K$ and
setting $\bB=(\bone,\cdots,\bone)_{p\times K}$. Further suppose we
know ideally that $\cov(\bff)=\var(f_1)I_K$. Then we have
\[ \Sig_n=I_{p}+\var(f_1)K\bone\bone\t\quad\text{and}\quad
\hSig_n=I_{p}+\hvar(f_1)K\bone\bone\t. \] It is easy to see that
$\lambda_1(\Sig_n)=\var(f_1)pK+1$, $\lambda_k(\Sig_n)=1$,
$k=2,\cdots,p$ and $\lambda_1(\hSig_n)=\hvar(f_1)pK+1$,
$\lambda_k(\hSig_n)=1$, $k=2,\cdots,p$. Thus,
\[ \max_{1\leq k\leq
p}\left|\lambda_k(\hSig_n)-\lambda_k(\Sig_n)\right|=\left|\hvar(f_1)-\var(f_1)\right|pK
=O_P(n^{-1/2}pK).
\]
Therefore, $\hSig$ here attains the optimal uniform weak convergence
rate of eigenvalues.

Theorem 1 shows that the factor structure does not give much
advantage in estimating $\Sig$.  The next theorem shows that when
$\Sig^{-1}$ is involved, the rate of convergence is improved.

\medskip

\begin{thm} (Rates of convergence under norm $\|\cdot\|_\Sig$).\quad\em Suppose that
$K=O(n^{\alpha_1})$ and $p=O(n^\alpha)$. Under conditions (A)--(C),
we have $\left\|\hSig-\Sig\right\|_\Sig=O_P(n^{-\beta/2})$ with
$\beta=\min\left(1-2\alpha_1,2-\alpha-\alpha_1\right)$ and
$\left\|\hSig\sam-\Sig\right\|_\Sig=O_P(n^{-\beta_1/2})$ with
$\beta_1=1-\max(\alpha,3\alpha_1/2,\\3\alpha_1-\alpha)$.
\end{thm}

\medskip

It is easy to show that $\beta > \beta_1$ whenever $\alpha >  2
\alpha_1$ and $\alpha_1<1$. Hence, the sample covariance matrix
$\hSig\sam$ has slower convergence. An interesting case is $K=O(1)$.
In this case, under the norm $\|\cdot\|_\Sig$, $\hSig$ has
convergence rate $n^{-\beta/2}$ with $\beta=\min(1,2-\alpha)$,
whereas $\hSig\sam$ has slower convergence rate $n^{-\beta_1/2}$
with $\beta_1 = 1 - \alpha$. In particular, when $\alpha \leq 1$,
$\hSig$ is root-$n$-consistent under $\|\cdot\|_\Sig$. This can be
shown to be optimal by some calculations using a specific factor
model mentioned above.

\bigskip

\begin{thm} (Rates of convergence of inverse under Frobenius norm).\quad\em Under conditions
(A)--(C), we have
\[ \left\|\hSig_n^{-1}-\Sig_n^{-1}\right\|=o_P\{(p^2K^4\log
n/n)^{1/2}\}, \] whereas
\[ \left\|\hSig\sam^{-1}-\Sig_n^{-1}\right\|=o_P\{(p^4K^2\log
n/n)^{1/2}\}. \]
\end{thm}

From this theorem, we see that when $K=o(p)$, $\hSig^{-1}$ performs
much better than $\hSig\sam^{-1}$. As expected, they perform roughly
the same in the extreme case where $K$ is proportional to $p$. It is
very pleasing that under an additional assumption (C), $\hSig^{-1}$
has a consistency rate slightly slower than $\hSig$ under the
Frobenius norm, since $\hSig^{-1}$ involves the inverse of the
$K\times K$ sample covariance matrix of $\bff$. The consistency
result of $\hSig\sam^{-1}$ is implied by that of $\hSig\sam$, thanks
to a simple inequality in matrix theory on inverses under
perturbation. However, the consistency result of $\hSig^{-1}$ needs
a very delicate analysis of inverse matrices. This theorem will be
used in Section 3.1 to examine the variance of a mean-variance
optimal portfolio.

\bigskip

Before going further, we first introduce some standard notation. Let
$\bA=(a_{ij})$ be a $q\times r$ matrix and denote by $\vecc(\bA)$
the $qr\times1$ vector formed by stacking the $r$ columns of $\bA$
underneath each other in the order from left to right. In
particular, for any $d\times d$ symmetric matrix $\bA$, we denote by
$\vech(\bA)$ the $d(d+1)/2\times1$ vector obtained from $\vecc(\bA)$
by removing the above-diagonal entries of $\bA$. It is not difficult
to see that there exists a unique $d^2\times d(d+1)/2$ matrix $D_d$
of zeros and ones such that
\[ D_d\ \vech(\bA)=\vecc(\bA) \]
for any $d\times d$ symmetric matrix $\bA$. $D_d$ is called the
duplication matrix of order $d$. Clearly, for any $d\times d$
symmetric matrix $\bA$, we have
\[ P_D\vecc(\bA)=\vech(\bA), \]
where $P_D=\left(D\t D\right)^{-1}D\t$. For any $q\times r$ matrix
$\bA_1=(a_{ij})$ and $s\times t$ matrix $\bA_2$, we define their
Kronecker product $\bA_1\otimes\bA_2$ as the $qs\times rt$ matrix
$(a_{ij}\bA_2)$.

\bigskip

\begin{thm} (Asymptotic normality).\quad\em Under conditions (A),
(B), and (D), if $p\rightarrow\infty$ as $n\rightarrow\infty$, then
the estimator $\hSig$ satisfies
\[ \sqrt{n}\ \vech\left[p^{-2}\bB_n\t\left(\hSig_n-\Sig_n\right)\bB_n\right]\toD\mathcal{N}\left(0,G\right), \]
where $G=P_D\left(\bA\otimes \bA\right)DHD\t\left(\bA\otimes
\bA\right)P_D\t$, $H=\cov\left[\vech\left(U\right)\right]$ with
$U=(u_{ij})_{K\times K}$ and \[
\cov\left(u_{ij},u_{kl}\right)=\kappa^{ijkl}+\kappa^{ik}\kappa^{jl}+
\kappa^{il}\kappa^{jk}, \] $\kappa^{i_1\cdots i_r}$ is the central
moment $E\left[(f_{i_1}-Ef_{i_1})\cdots(f_{i_r}-Ef_{i_r})\right]$ of
$\bff=(f_1,\cdots,f_K)\t$, $D$ is the duplication matrix of order
$K$, and $P_D=\left(D\t D\right)^{-1}D\t$.
\end{thm}

\medskip

When $\bff$ has a $K$-variate normal distribution with covariance
matrix $(\sigma_{ij})_{K\times K}$, the matrix $H$ in Theorem 4 is
determined by
\[ \cov\left(u_{ij},u_{kl}\right)=\sigma_{ik}\sigma_{jl}+\sigma_{il}\sigma_{jk}. \]

The diverging dimensionality takes care of a trouble term in
establishing asymptotic normality. However, in the finite
dimensional setting, one can only show asymptotic normality when
$\bff$ has mean $\bzero$, where $\cov(\bff)$ can be estimated as
$\hcov(\bff)=n^{-1}\bX\bX\t$, and in general, $\hSig$ may have no
asymptotic normality because the term
$\bX\bone\bone\t\bX\t\left(\bX\bX\t\right)^{-1}\bX$ may not have a
limiting behavior as $n\rightarrow\infty$ (at least it is not clear
now). This is an interesting phenomenon in the presence of diverging
dimensionality.

\bigskip

\setcounter{equation}{0} \setcounter{subsect}{0}
\setcounter{subsubsect}{0}

\sect {\bf Impacts on portfolio allocation and risk
management.}\quad In this section we examine the impacts of
covariance matrix estimation on portfolio allocation and risk
management, respectively.

\bigskip

\subsect {\em Impact on portfolio allocation}.\quad For practical
use in portfolio allocation, one would expect that the optimal
portfolio constructed from the covariance matrix estimated from the
history should not deviate too much from the true one. So we examine
the behavior of the optimal portfolio constructed using $\hSig$
estimated from historical data.

Markowitz (1952) defines the mean-variance optimal portfolio as the
solution $\bxi_n\in\mathbb{R}^{p}$ to the following minimization
problem
\begin{align} \label{036}
&\min_{\bxi}\ \bxi\t\Sig_n\bxi\\
\nonumber &\text{Subject to }\bxi\t\bone=1\text{ and
}\bxi\t\bmu_n=\gamma_n,
\end{align}
where $\bone$ is a $p\times1$ vector of ones,
$\bmu_n=E\left(\by\right)$, and $\gamma_n$ is the expected rate of
return imposed on the portfolio. It is well known that Markowitz's
optimal portfolio [see Markowitz (1959), Cochrane (2001), or
Campbell, Lo and MacKinlay (1997)] is
\begin{equation} \label{037}
\bxi_n=\frac{\phi_n-\gamma_n
\psi_n}{\varphi_n\phi_n-\psi_n^2}\Sig_n^{-1}\bone+\frac{\gamma_n
\varphi_n-\psi_n}{\varphi_n\phi_n-\psi_n^2}\Sig_n^{-1}\bmu_n
\end{equation}
with $\varphi_n=\bone\t\Sig_n^{-1}\bone$,
$\psi_n=\bone\t\Sig_n^{-1}\bmu_n$, and
$\phi_n=\bmu_n\t\Sig_n^{-1}\bmu_n$, and its variance is
\begin{equation} \label{038}
\bxi_n\t\Sig_n\bxi_n=\frac{\varphi_n\gamma_n^2-2\psi_n\gamma_n+\phi_n}{\varphi_n\phi_n-\psi_n^2}.
\end{equation}
Denote by $\bxi_{ng}$ the $\bxi_n$ in (\ref{037}) with $\gamma_n$
replaced by $\psi_n/\varphi_n$. The global minimum variance without
constraint on the expected return is
\begin{equation} \label{137}
\bxi_{ng}\t\Sig_n\bxi_{ng}=\varphi_n^{-1},
\end{equation}
which is attained in (\ref{038}) when $\gamma_n=\psi_n/\varphi_n$.

Based on the history, we can construct $\hSig_n$ as before. Also, we
have a substitution estimator
$\hmu_n=\hB_nn^{-1}(\bff_1+\cdots+\bff_n)$ of the mean vector
$\bmu_n$. As above, we can define estimators $\hxi_n$, $\hxi_{ng}$
and $\widehat{\varphi}_n$, $\widehat{\psi}_n$, $\widehat{\phi}_n$
with $\Sig_n$ and $\bmu_n$ replaced by $\hSig_n$ and $\hmu_n$,
respectively.

It is interesting to study the deviation of the constructed optimal
portfolio $\hxi_n$ and the globally optimal portfolio $\hxi_{ng}$
from the theoretical ones, say, $\bxi_n$ and $\bxi_{ng}$. But here
we do not pursue in this direction because it is more valuable to
study the risk associated with them. Therefore, we only examine the
behavior of the minimum variance $\hxi_n\t\hSig_n\hxi_n$ and global
minimum variance $\hxi_{ng}\t\hSig_n\hxi_{ng}$ in this section.

\bigskip

\begin{thm} (Weak convergence of global minimum variance).\quad\em Suppose that
all the $\varphi_n$'s are bounded away from zero. Under conditions
(A)--(C), we have
\[
\hxi_{ng}\t\hSig_n\hxi_{ng}-\bxi_{ng}\t\Sig_n\bxi_{ng}=o_P\{(p^4K^4\log
n/n)^{1/2}\},
\]
whereas
\[
\hxi_{ng}\t\hSig\sam\hxi_{ng}-\bxi_{ng}\t\Sig_n\bxi_{ng}=o_P\{(p^6K^2\log
n/n)^{1/2}\}.
\medskip
\]
\end{thm}
\begin{thm} (Weak convergence to optimal portfolio).\quad\em Suppose that
$\varphi_n\phi_n-\psi_n^2$ are bounded away from zero and
$\varphi_n/(\varphi_n\phi_n-\psi_n^2)$,
$\psi_n/(\varphi_n\phi_n-\psi_n^2)$,
$\phi_n/(\varphi_n\phi_n-\psi_n^2)$, $\gamma_n$ are bounded. Under
conditions (A)--(C), we have
\[
\hxi_n\t\hSig_n\hxi_n-\bxi_n\t\Sig_n\bxi_n=o_P\{(p^4K^4\log
n/n)^{1/2}\}, \] whereas
\[
\hxi_n\t\hSig\sam\hxi_n-\bxi_n\t\Sig_n\bxi_n=o_P\{(p^6K^2\log
n/n)^{1/2}\}. \]
\end{thm}
\quad The assumptions on $\varphi_n$, $\psi_n$ and $\phi_n$ in
Theorems 5 and 6 are technical and reasonable. In view of
(\ref{137}), the assumption on $\varphi_n$ in Theorem 5 amounts to
saying that the global minimum variances are bounded across $n$. The
additional  assumptions in Theorem 6 can be understood in a similar
way in light of (\ref{038}). From the above two theorems, we see
that when $K=o(p)$, $\hSig$ performs much better than $\hSig\sam$
from the point of view of portfolio allocation. On the other hand,
we also see that dimensionality as well as number of factors can
only grow slowly with sample size so that the globally optimal
portfolio and the mean-variance optimal portfolio constructed using
estimated covariance matrix $\hSig$ or $\hSig\sam$ behave similarly
to theoretical ones. So high dimensionality does impose a great
challenge on portfolio allocation.

Our study reveals that for a large number of stocks, additional
structures are needed. For example, we may group assets according to
sectors and assume that the sector correlations are weak and
negligible. Hence, the covariance structure is block diagonal. Our
factor model approach can be used to estimate the covariance matrix
within a block, and our results continue to apply.

\bigskip

\subsect {\em Impact on risk management}.\quad Risk management is a
different story from portfolio allocation. As mentioned in Section
1.1, the smallest and largest eigenvalues of the covariance matrix
are related to the minimum and maximum variances of the selected
portfolio, respectively. Throughout this section, we fix a sequence
of selected portfolios $\bxi_n\in\mathbb{R}^{p}$ with
$\bxi_n\t\bone=1$ and $\bxi_n=O(1)\bone$. Here we impose the
condition $\bxi_n=O(1)\bone$ to avoid extreme short positions --
that is, some large negative components in $\bxi_n$. Then, the
variance of portfolio $\bxi_n$ is
\[ \var(\bxi_n\t\by)=\bxi_n\t\cov(\by)\bxi_n=\bxi_n\t\Sig_n\bxi_n. \]
The estimated risk associated with portfolio $\bxi_n$ is
$\bxi_n\t\hSig_n\bxi_n$. For practical use in risk management, we
need to examine the behavior of portfolio variance based on
$\hSig_n$ estimated from historical data.

\bigskip

\begin{thm} (Weak convergence of variance).\quad\em Under conditions (A) and (B), we have
\[ \bxi_n\t\hSig_n\bxi_n-\bxi_n\t\Sig_n\bxi_n=o_P\{(p^4K^2\log
n/n)^{1/2}\} \] and
\[ \bxi_n\t\hSig\sam\bxi_n-\bxi_n\t\Sig_n\bxi_n=o_P\{(p^4K^2\log
n/n)^{1/2}\}. \] On the other hand, if the portfolios $\bxi_n$'s
have no short positions, then we have
\[ \bxi_n\t\hSig_n\bxi_n-\bxi_n\t\Sig_n\bxi_n=o_P\{(p^2K^2\log
n/n)^{1/2}\} \] and
\[ \bxi_n\t\hSig\sam\bxi_n-\bxi_n\t\Sig_n\bxi_n=o_P\{(p^2K^2\log
n/n)^{1/2}\}.\medskip \]
\end{thm}
\quad From this theorem, we see that $\hSig$ behaves roughly the
same as the sample covariance matrix estimator $\hSig\sam$ in risk
management. This is essential for both covariance matrix estimators,
since risk management does not involve inverse of the covariance
matrix, but the covariance matrix itself. The above theorem is
implied by consistency results of $\hSig$ and $\hSig\sam$ under the
Frobenius norm in Theorem 1.

\bigskip
\setcounter{equation}{0} \setcounter{subsect}{0}
\setcounter{subsubsect}{0}

\sect {\bf A simulation study.}\quad In this section we use a
simulation study to illustrate and augment our theoretical results
and to verify finite-sample performance of the estimator $\hSig$ as
well as $\hSig^{-1}$. To this end, we fix sample size $n=756$, which
is the practical sample size of three-year daily financial data, and
we let dimensionality $p$ grow from low to high and ultimately
exceed sample size. As mentioned before, our primary concern is a
theoretical understanding of factor models with a diverging number
of variables and factors for the purpose of covariance matrix
estimation, but not comparison with other popular estimators. So we
compare performance of the estimator $\hSig$ only to that of sample
covariance matrix $\hSig\sam$. To contrast with $\hSig\sam$, we
examine the covariance matrix estimation errors of $\hSig$ and
$\hSig\sam$ under the Frobenius norm, the norm $\|\cdot\|_\Sig$
introduced in Section 2, and the Stein (or entropy) loss function
\[
L(\hSig,\Sig)=\tr\left(\hSig\Sig^{-1}\right)-\log\left|\hSig\Sig^{-1}\right|-p,
\]
which was proposed by James and Stein (1961). Meanwhile, we compare
estimation errors of $\hSig^{-1}$ and $\hSig\sam^{-1}$ under the
Frobenius norm. Furthermore, we evaluate estimated variances of
optimal portfolios with expected rate of return $\gamma_n=10\%$
based on $\hSig$ and $\hSig\sam$ by comparing their mean-squared
errors (MSEs). For the estimated global minimum variances, we also
compare their MSEs. Moveover, we examine MSEs of estimated variances
of the equally weighted portfolio $\bxi_p=(1/p,\cdots,1/p)$, based
on $\hSig$ and $\hSig\sam$, respectively.

For simplicity, we fix $K=3$ in our simulation and consider the
three-factor model
\begin{equation} \label{145}
Y_{pi}=b_{pi1}f_1+b_{pi2}f_2+b_{pi3}f_3+\veps_i,\quad i=1,\cdots,p.
\end{equation}
Here, we use the first subscript $p$ to stress that the three-factor
model varies across dimensionality $p$. As before, we let
$\by=(Y_1,\cdots,Y_p)\t$ and $\bff=(f_1,f_2,f_3)\t$. The Fama-French
three-factor model [Fama and French (1993)] is a practical example
of model (\ref{145}). To make our simulation more realistic, we take
the parameters from a fit of the Fama-French three-factor model.

In the Fama-French three-factor model, $Y_i$ is the excess return of
the $i$-th stock or portfolio, $i=1,\cdots,p$. The first factor
$f_1$ is the excess return of the proxy of the market portfolio,
which is the value-weighted return on all NYSE, AMEX and NASDAQ
stocks (from CRSP) minus the one-month Treasury bill rate (from
Ibbotson Associates). The other two factors are constructed using
six value-weighted portfolios formed on size and book-to-market.
Specifically, the second factor $f_2$, \verb=SMB= (Small Minus Big),
\begin{align*}
\text{SMB}&=1/3\left(\text{Small Value}+\text{Small
Neutral}+\text{Small Growth}\right)\\
&\quad-1/3\left(\text{Big Value}+\text{Big Neutral}+\text{Big
Growth}\right)
\end{align*}
is the average return on the three small portfolios minus the
average return on the three big portfolios, and the third factor
$f_3$, \verb=HML= (High Minus Low),
\begin{align*}
\text{HML}&=1/2\left(\text{Small Value}+\text{Big Value}\right)\\
&\quad-1/2\left(\text{Small Growth}+\text{Big Growth}\right)
\end{align*}
is the average return on the two value portfolios minus the average
return on the two growth portfolios. See their website
\verb=http://mba.tuck.dartmouth.edu/pages/faculty=\\
\verb=/ken.french/data_library.html= for more details about their
three factors and the data sets of the three factors, risk free
interest rates, and returns of many constructed portfolios.

We first fit three-factor model (\ref{145}) with $n=756$ and $p=30$
using the three-year daily data of 30 Industry Portfolios from May
1, 2002 to Aug. 29, 2005, which are available at the above website.
Then, as in (\ref{105}), we get 30 estimated factor loading vectors
$\hb_1=(b_{11},b_{12},b_{13}),\cdots,\hb_{30}=(b_{30,1},b_{30,2},b_{30,3})$
and 30 estimated standard deviations $\hsig_1,\cdots,\hsig_{30}$ of
the errors, where $\hb_i$ and $\hsig_i$ correspond to the $i$-th
portfolio, $i=1,\cdots,30$. The sample average of
$\hsig_1,\cdots,\hsig_{30}$ is 0.66081 with a sample standard
deviation 0.3275. We report in Table 1 the sample means and sample
covariance matrices of $\bff$ and $\hb$ denoted by $\mu_{\bff}$,
$\mu_{\bb}$ and $\cov_{\bff}$, $\cov_{\bb}$, respectively.

\vspace{0.1 in}
\begin{table}[tph]
\begin{center}
{\small \textsc{Table 1}\\
\textit{Sample means and sample covariance matrices of $\bff$ and
$\hb$}\vspace{0.1 in}}
\begin{tabular}
[c]{c|lll}\hline $\mu_{\bff}$ & & $\cov_{\bff}$ & \\
\hline 0.023558 & \ \ 1.2507 & -0.034999 & -0.20419 \\
0.012989 & \ -0.034999 & \ 0.31564 & -0.0022526 \\
0.020714 & \ -0.20419 & -0.0022526 & \ 0.19303 \\
\hline\hline $\mu_{\bb}$ & & $\cov_{\bb}$ & \\
\hline 0.78282 & \ 0.029145 & \ 0.023873 & \ 0.010184 \\
0.51803 & \ 0.023873 & \ 0.053951 & -0.006967 \\
0.41003 & \ 0.010184 & -0.006967 & \ 0.086856 \\\hline
\end{tabular}
\end{center}
\end{table}

For each simulation, we carry out the following steps:
\begin{itemize}
\item We first generate a random sample of $\bff=(f_1,f_2,f_3)\t$ with
size $n=756$ from the trivariate normal distribution
$\mathcal{N}\left(\mu_{\bff},\cov_{\bff}\right)$.

\item Then, for each dimensionality $p$ increasing from 16 to 1000 with
increment 20, we do the following.

\item Generate $p$ factor loading vectors
$\bb_1,\cdots,\bb_p$ as a random sample of size $p$ from the
trivariate normal distribution
$\mathcal{N}\left(\mu_{\bb},\cov_{\bb}\right)$.

\item Generate $p$
standard deviations $\sigma_1,\cdots,\sigma_p$ of the errors as a
random sample of size $p$ from a gamma distribution
$G(\alpha,\beta)$ conditional on being bounded below by a threshold
value. The threshold for the standard deviations of errors is
required in accordance with condition (C) in Section 2.1, and it is
set to 0.1950 in our simulation because we find $\min_{1\leq i\leq
30}\hsig_i=0.1950$. Note that $G(\alpha,\beta)$ has mean
$\alpha\beta$ and standard deviation $\alpha^{1/2}\beta$, and its
conditional mean and conditional second moment on falling above
0.1950 can be approximated respectively by \[
\left(\alpha\beta-\frac{0.1950}{2}p\right)/\left(1-p\right)\quad
\text{and} \quad
\left(\alpha\beta^2+\alpha^2\beta^2-\frac{0.1950^2}{2}p\right)/\left(1-p\right),\]
where $p$ is the probability of falling below 0.1950 under
$G(\alpha,\beta)$. By matching the mean 0.66081 and standard
deviation 0.3275 for $G(\alpha_0,\beta_0)$, we obtain
$\alpha_0=4.0713$ and $\beta_0=0.1623$. Therefore, following the
above approximations, by recursively matching the conditional mean
0.66081 and conditional second moment $0.3275^2+0.66081^2=0.54393$
for $G(\alpha,\beta)$, we finally get $\alpha=3.3586$ and
$\beta=0.1876$.

\item After getting $p$
standard deviations $\sigma_1,\cdots,\sigma_p$ of the errors, we
generate a random sample of $\bveps=(\veps_1,\cdots,\veps_p)\t$ with
size $n=756$ from the $p$-variate normal distribution
$\mathcal{N}\left(0,\diag\left(\sigma_1^2,\cdots,\sigma_p^2\right)\right)$.

\item Then from model (\ref{145}), we get a random sample of
$\by=(Y_1,\cdots,Y_p)\t$ with size $n=756$.

\item Finally, we compute estimated covariance matrices
$\hSig$ and $\hSig\sam$, as well as $\hSig^{-1}$ and
$\hSig\sam^{-1}$, and record the errors in the aforementioned
measures. Meanwhile, we calculate MSEs of estimated variances of the
optimal portfolios with $\gamma_n=10\%$ as well as MSEs of estimated
global minimum variances based on $\hSig$ and $\hSig\sam$,
respectively. Also, we record MSEs of estimated variances of the
equally weighted portfolio based on $\hSig$ and $\hSig\sam$,
respectively.
\end{itemize}
We repeat the above simulation 500 times and report the mean-square
errors as well as the standard deviations of those errors.

\begin{singlespace}
\begin{figure} \centering
\begin{center}%
\begin{tabular}
[c]{cc}%
{\includegraphics[ trim=-0.000000in 0.000000in 0.000000in
-0.189256in, height=2in, width=2.6in
]%
{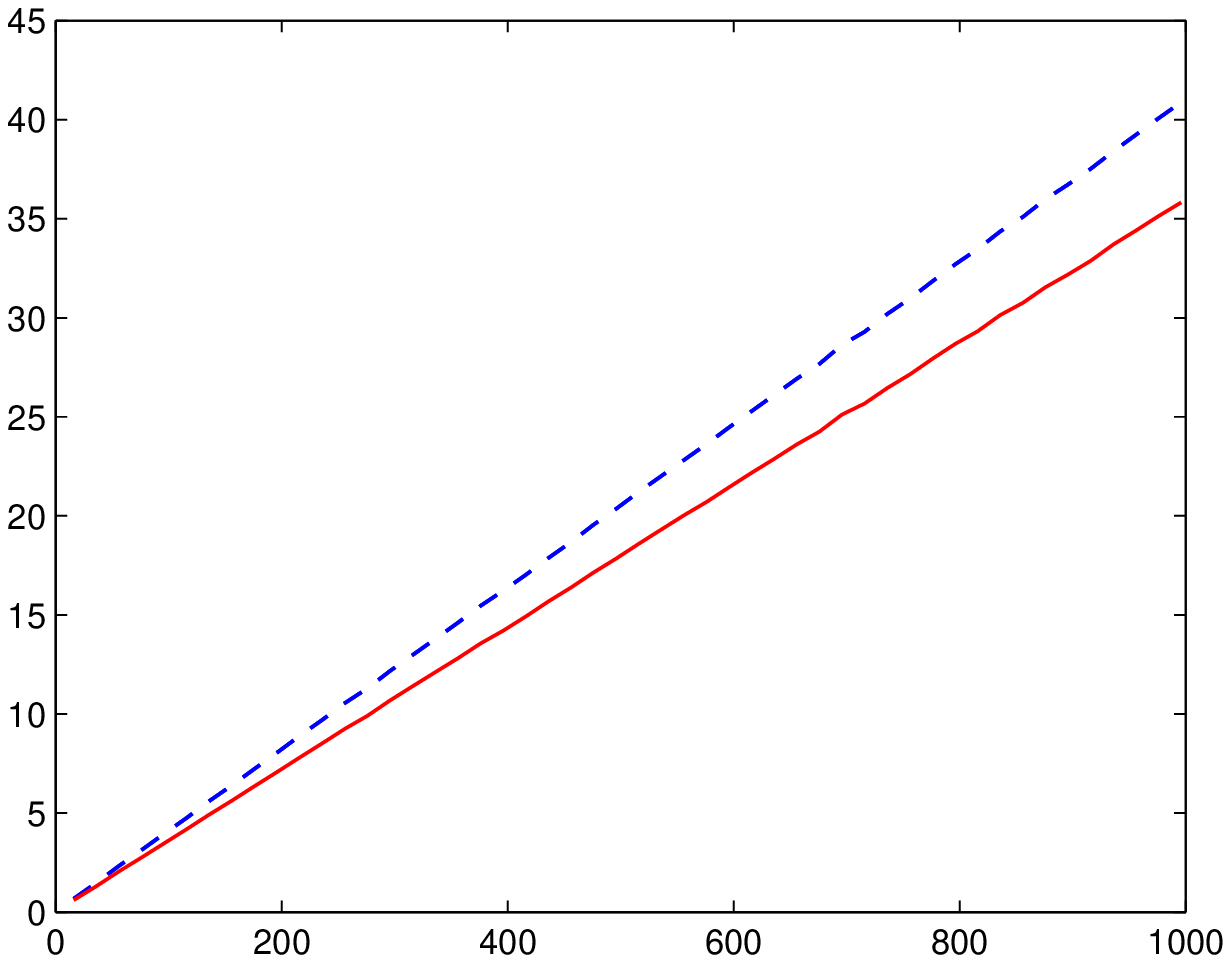}%
}%
& {\includegraphics[ trim=0.000000in 0.000000in 0.000000in
-0.189256in, height=2in, width=2.6in
]%
{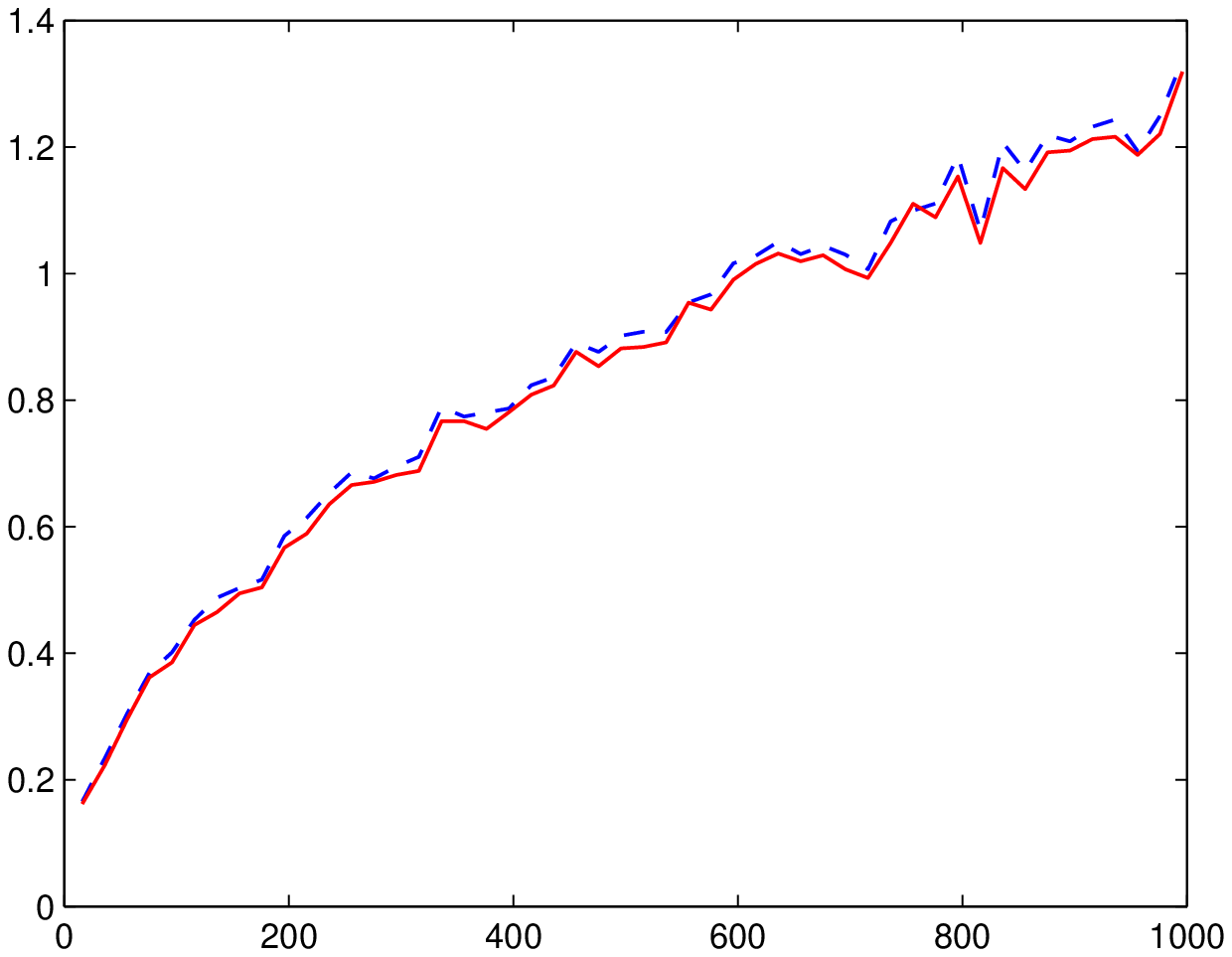}%
}%
\\
\scriptsize (a) & \scriptsize (b)\\
{\includegraphics[ trim=0.000000in 0.000000in 0.000000in
-0.189256in, height=2in, width=2.6in
]%
{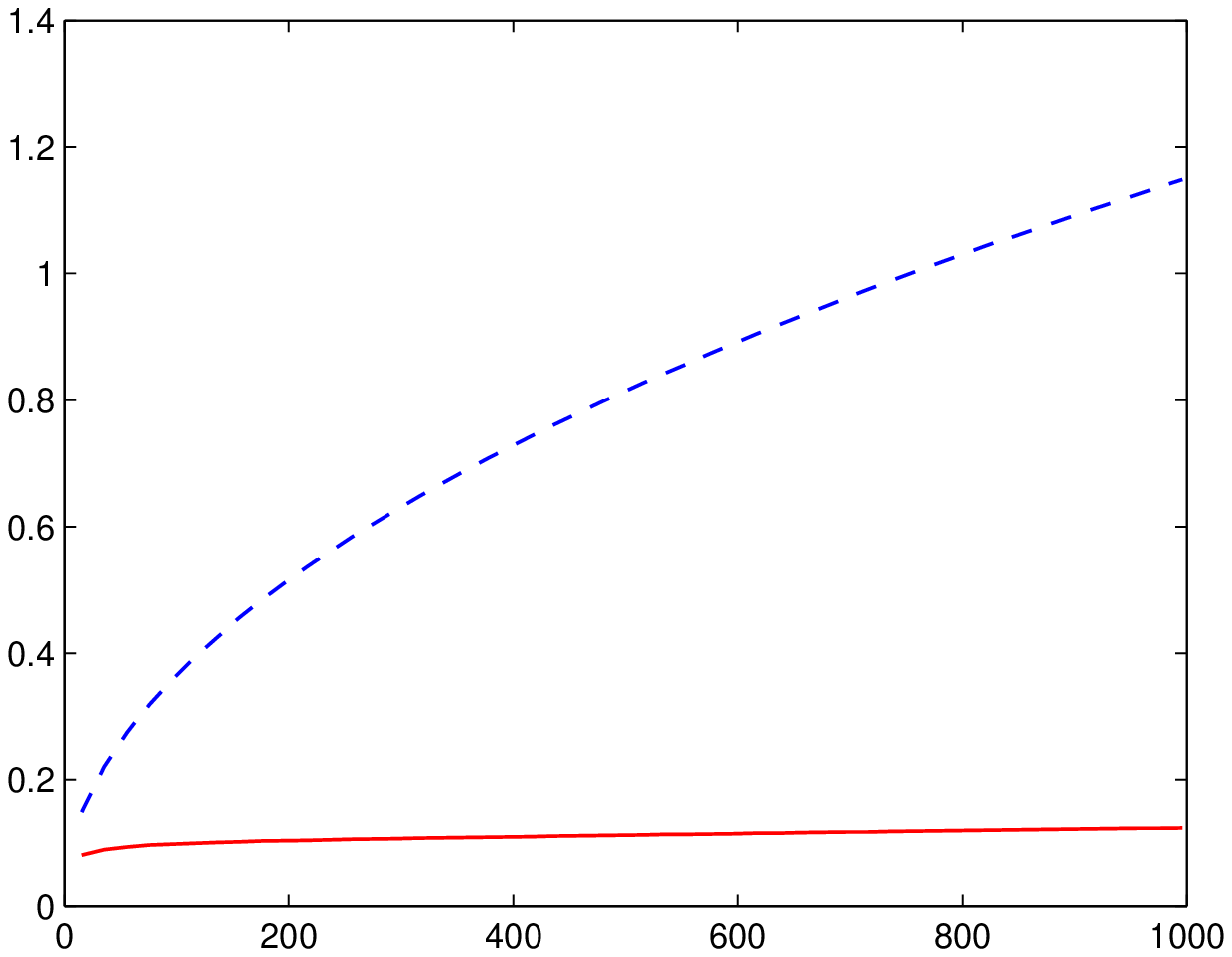}%
}%
& {\includegraphics[ trim=0.000000in 0.000000in 0.000000in
-0.189256in, height=2in, width=2.6in
]%
{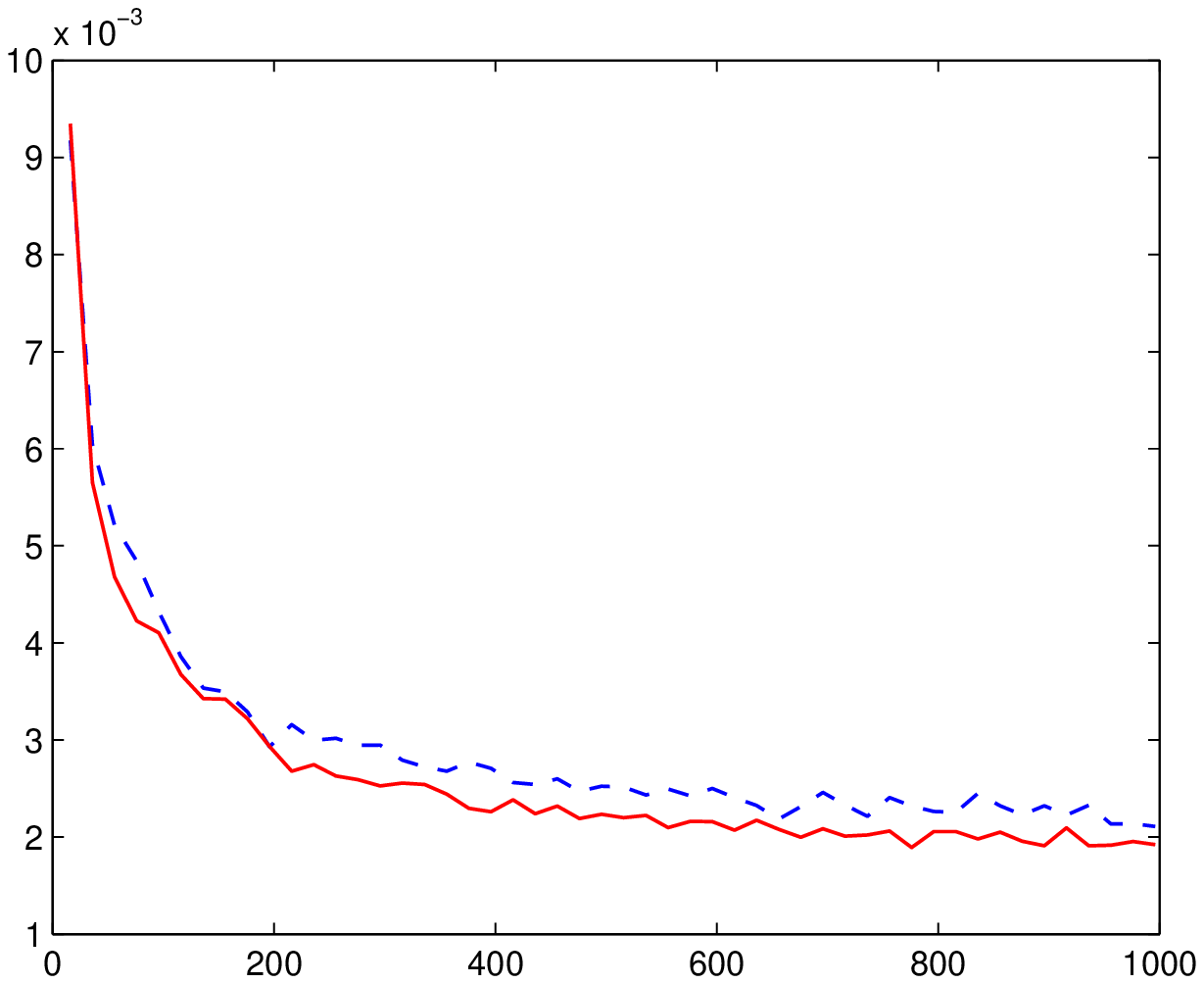}%
}%
\\
\scriptsize (c) & \scriptsize (d)\\
{\includegraphics[ trim=0.000000in 0.000000in 0.000000in
-0.189256in, height=2in, width=2.6in
]%
{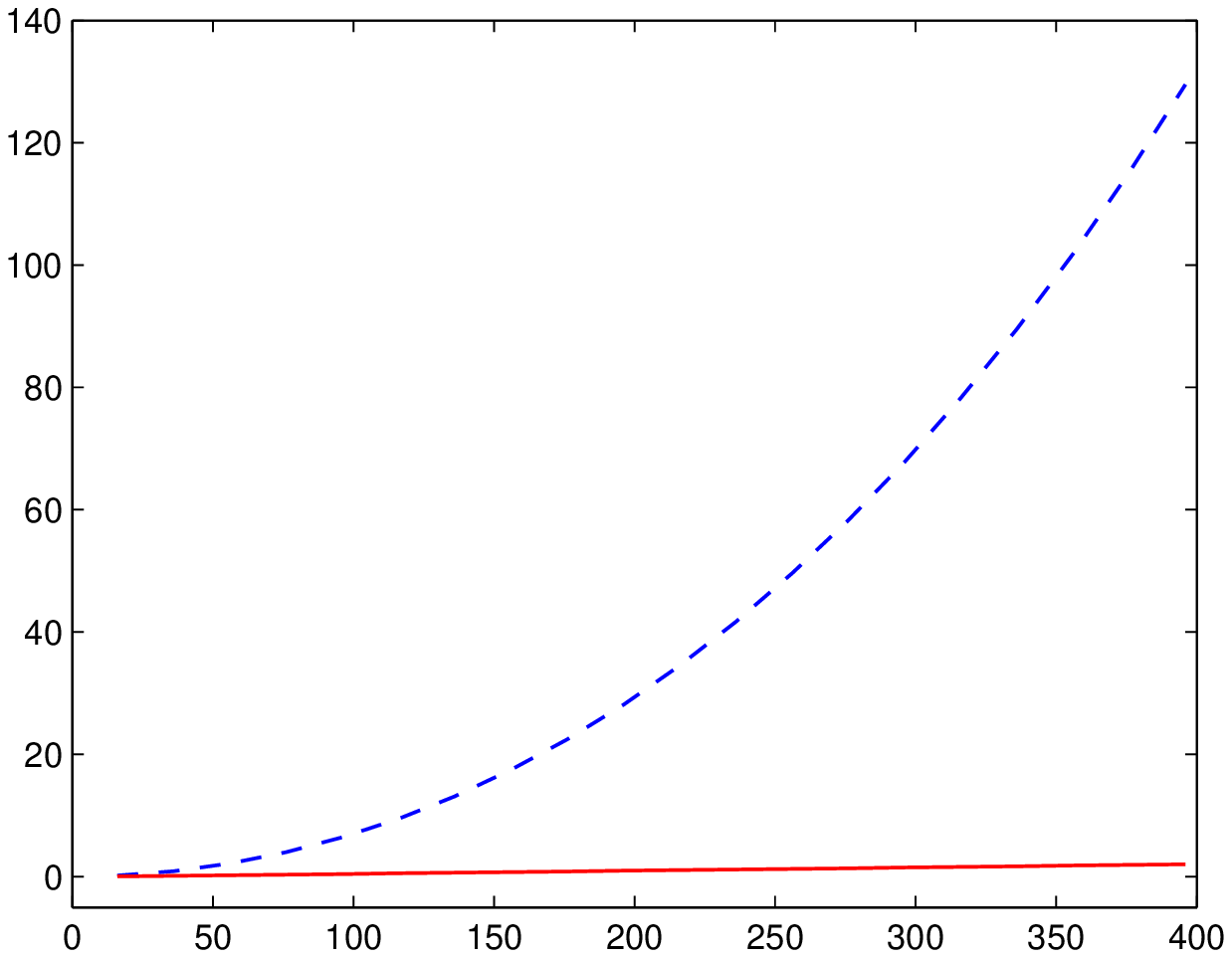}%
}%
& {\includegraphics[ trim=0.000000in 0.000000in 0.000000in
-0.189256in, height=2in, width=2.6in
]%
{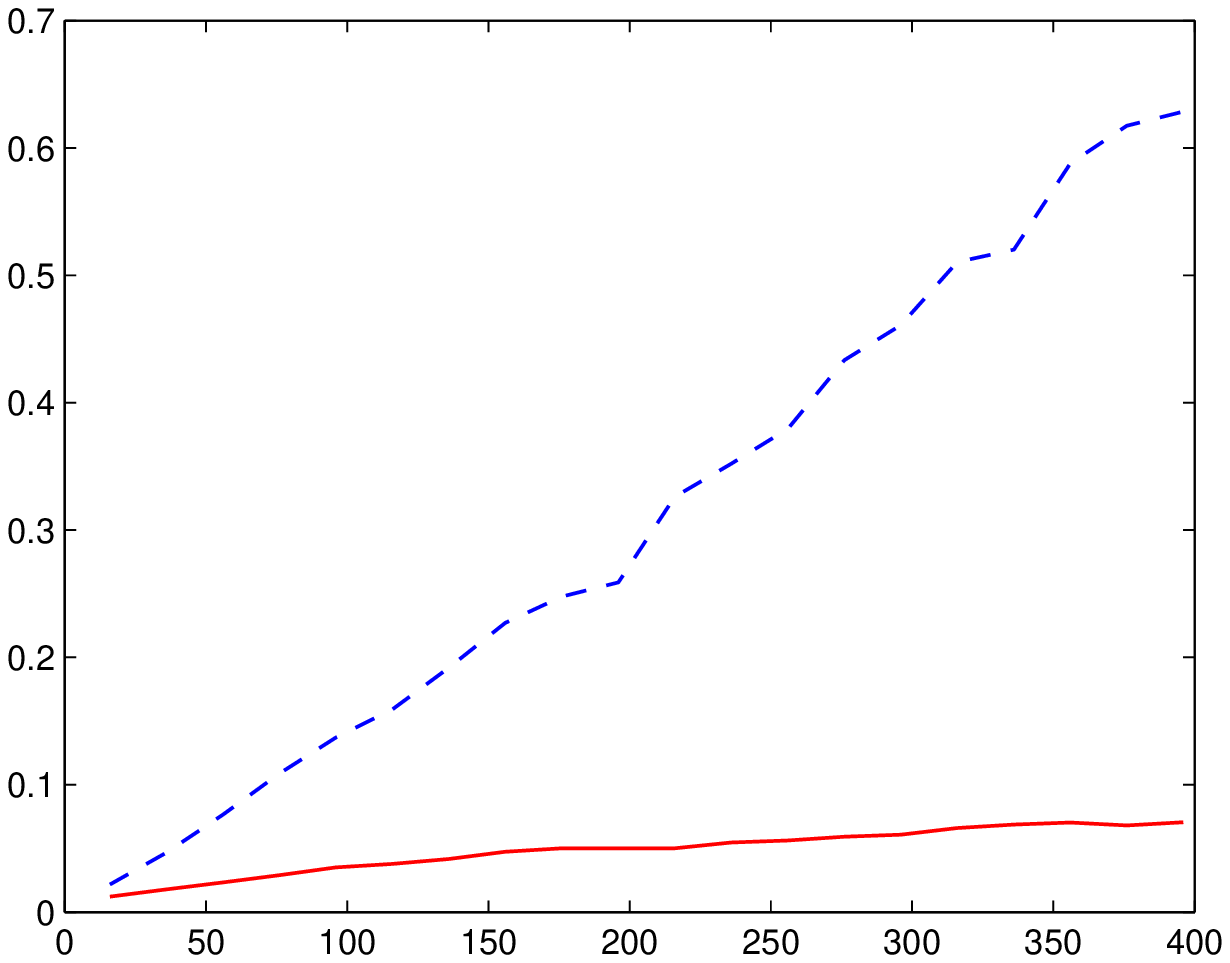}%
}%
\\
\scriptsize (e) & \scriptsize (f)
\end{tabular}%
\caption{{\scriptsize (a), (c) and (e): The averages of errors over
500 simulations for $\hSig$ (solid curve) and $\hSig\sam$ (dashed
curve) against $p$ under Frobenius norm, norm $\|\cdot\|_\Sig$ and
entropy losses, respectively. (b), (d) and (f): Corresponding
standard deviations of errors over 500 simulations for $\hSig$
(solid curve) and $\hSig\sam$ (dashed curve).}}%
\end{center}%
\end{figure}%
\end{singlespace}

\begin{singlespace}
\begin{figure} \centering
\begin{center}%
\begin{tabular}
[c]{cc}%
{\includegraphics[ trim=-0.000000in 0.000000in 0.000000in
-0.189256in, height=2in, width=2.6in
]%
{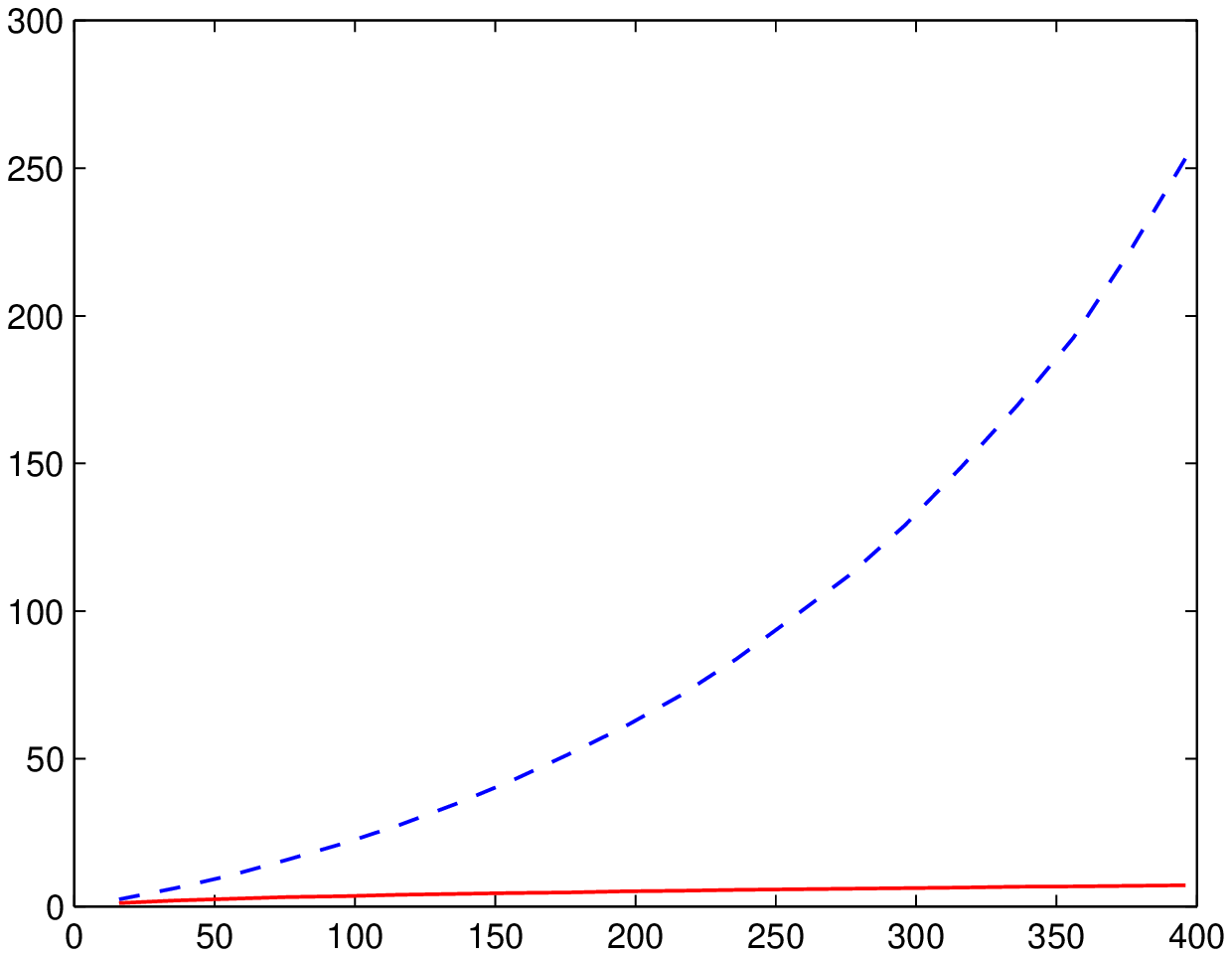}%
}%
& {\includegraphics[ trim=0.000000in 0.000000in 0.000000in
-0.189256in, height=2in, width=2.6in
]%
{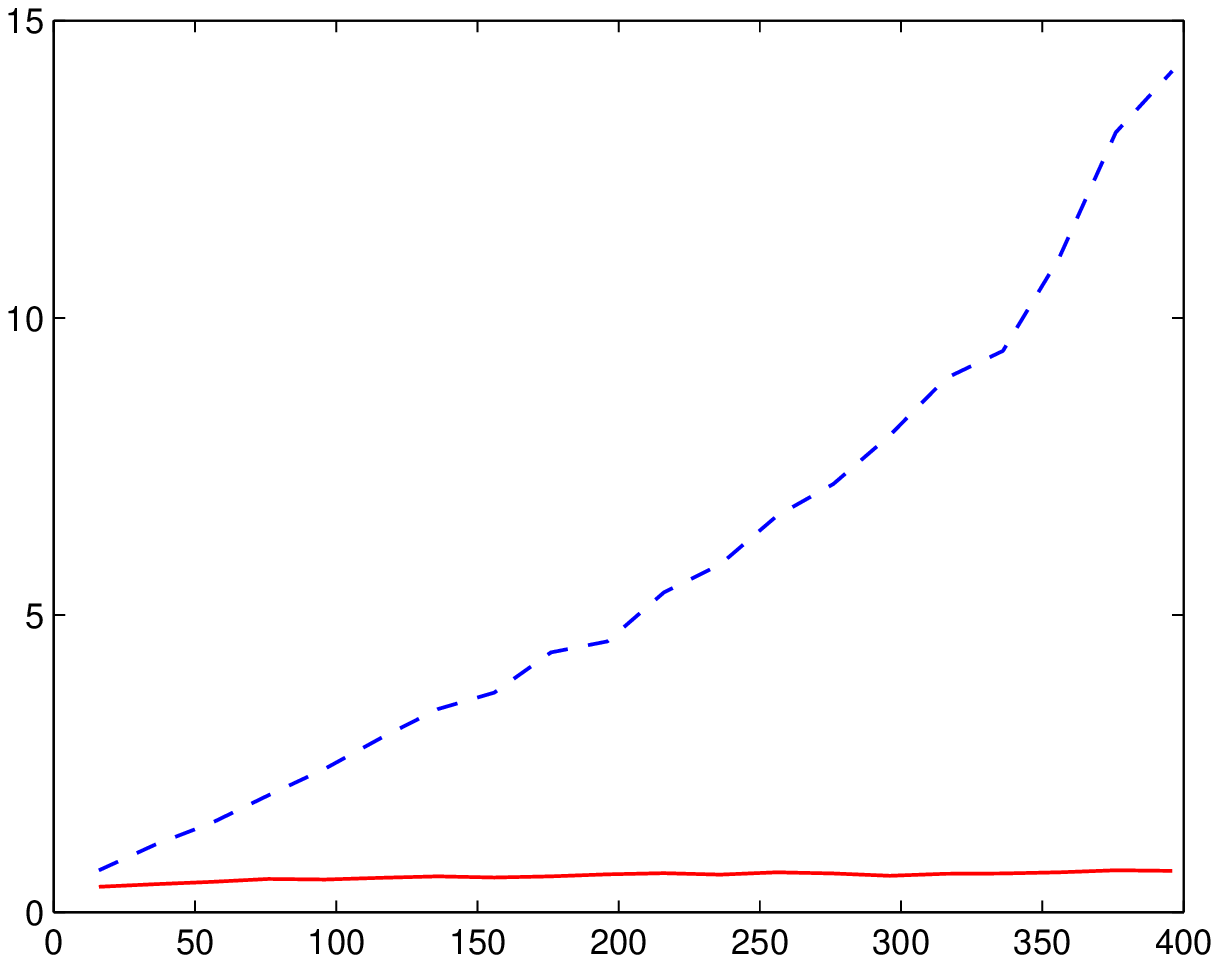}%
}%
\\
\scriptsize (a) & \scriptsize (b)
\end{tabular}%
\caption{{\scriptsize (a) The averages of errors under Frobenius
norm over 500 simulations for $\hSig^{-1}$ (solid curve) and
$\hSig\sam^{-1}$ (dashed curve) against $p$. (b) Corresponding
standard deviations of errors under Frobenius norm.}}%
\end{center}%
\end{figure}%
\end{singlespace}

In Figures 1--4, solid curves and dashed curves correspond to
$\hSig$ and $\hSig\sam$, respectively. Figure 1 presents the
averages and the standard deviations of their estimation errors
under the Frobenius norm, norm $\|\cdot\|_\Sig$, and entropy loss
against dimensionality $p$, respectively. Figure 2 depicts the
averages and the standard deviations of estimation errors of
$\hSig^{-1}$ and $\hSig\sam^{-1}$ under the Frobenius norm against
$p$. We report in Figure 3 MSEs of estimated variances of the
optimal portfolios with $\gamma_n=10\%$ as well as MSEs of estimated
global minimum variances using $\hSig$ and $\hSig\sam$ against $p$.
Figure 4 presents MSEs of estimated variances of the equally
weighted portfolio based on $\hSig$ and $\hSig\sam$ against $p$.

Recall that both the sample size $n$ and the number of factors $K$
are kept fixed across $p$ in our simulation. From Figures 1--4, we
observe the following:
\begin{itemize}
\item By comparing corresponding averages and standard
deviations of the errors shown in Figures 1 and 2, we see that the
Monte-Carlo errors are negligible.

\item Figure 1(a) shows that under the
Frobenius norm, $\hSig$ performs roughly the same as (slightly
better than) $\hSig\sam$, which is consistent with the results in
Theorem 1.  Nevertheless, this is a surprise and is against the
conventional wisdom.

\item Figure 1(c) reveals that under norm
$\|\cdot\|_\Sig$, $\hSig$ performs much better than $\hSig\sam$,
which is consistent with the results in Theorem 2. In particular, we
see that the estimation errors of $\hSig$ under norm
$\|\cdot\|_\Sig$ are roughly at the same level across $p$. Recall
that sample size $n$ is fixed as 756 here. Thus, this is in line
with the root-$n$-consistency of $\hSig$ under norm $\|\cdot\|_\Sig$
when $p=O(n)$ shown in Theorem 2. Also, the apparent growth pattern
of estimation errors in $\hSig\sam$ with $p$ is in accordance with
its $(n/p)^{1/2}$-consistency under norm $\|\cdot\|_\Sig$ shown in
Theorem 2.

\item Figure 1(e) shows that under entropy loss,
$\hSig$ significantly outperforms $\hSig\sam$, which strongly
supports the factor-model based estimator $\hSig$ over the sample
one $\hSig\sam$. We only report the results for $p$ truncated at
400. This is because for larger $p$, sample covariance matrices
$\hSig\sam$ are nearly singular with a big chance in the simulation,
which results in extremely large entropy losses.

\item From Figure 2(a), we see that under the
Frobenius norm, the estimator $\hSig^{-1}$ significantly outperforms
$\hSig\sam^{-1}$, which is in line with the results in Theorem 3.

\item Figures 3(a) and 3(b) demonstrate convincingly that $\hSig$
outperforms $\hSig\sam$ in portfolio allocation. These results are
in accordance with Theorems 5 and 6. One may notice that in Figure
3(a), the MSEs are relatively large in magnitude for small $p$ and
then tend to stabilize when $p$ grows large. This is because in our
settings for the simulation, for small $p$ the term
$\varphi_n\phi_n-\psi_n^2$ is relatively small compared to
$\varphi_n\gamma_n^2-2\psi_n\gamma_n+\phi_n$, which results in large
variance of the optimal portfolio. The behavior of the MSEs for
large $p$ is essentially due to self-averaging in the
dimensionality. Figures 3(b) can be interpreted in the same way.

\item  Figure 4 reveals that the factor-model based approach and the
sample approach have almost the same performance in risk management,
which is consistent with Theorem 7. The high-dimensionality behavior
is essentially due to self-averaging as in Figure 3(a).
\end{itemize}

\begin{singlespace}
\begin{figure} \centering
\begin{center}%
\begin{tabular}
[c]{cc}%
{\includegraphics[ trim=0.000000in 0.000000in 0.000000in
-0.189256in, height=2in, width=2.6in
]%
{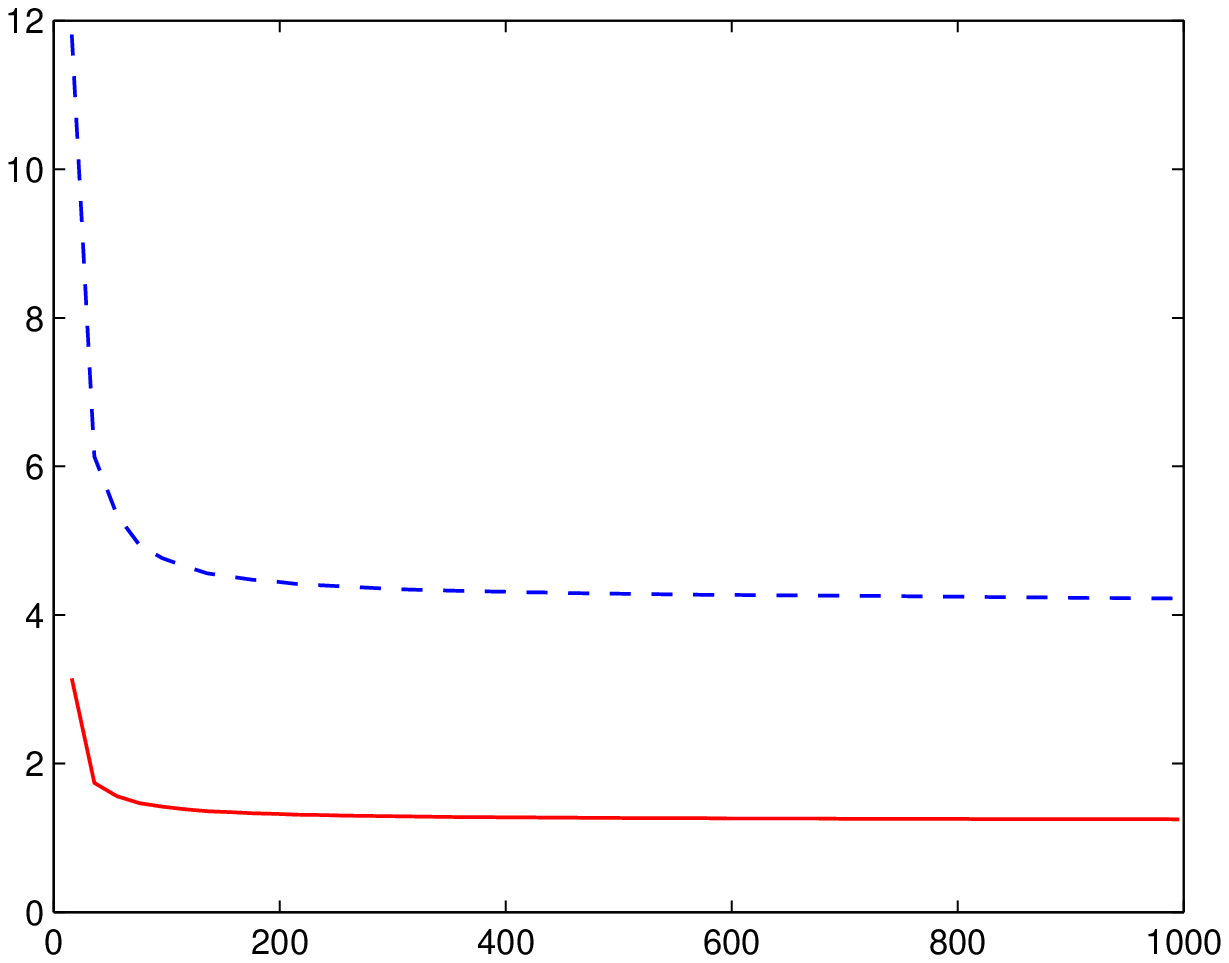}%
}%
& {\includegraphics[ trim=0.000000in 0.000000in 0.000000in
-0.189256in, height=2in, width=2.6in
]%
{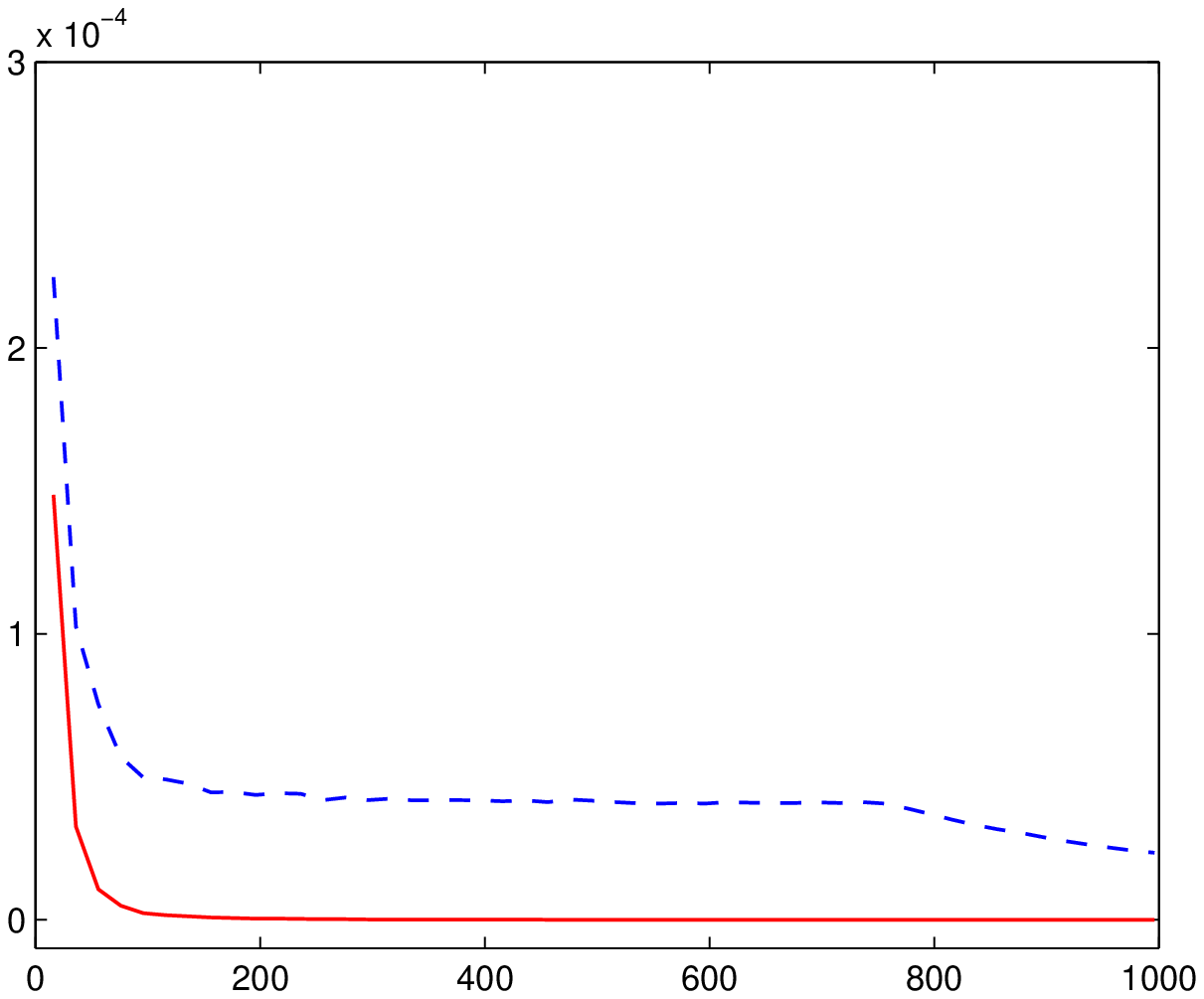}%
}%
\\
\scriptsize (a) & \scriptsize (b)
\end{tabular}%
\caption{{\scriptsize (a) The MSEs of estimated variances of the
optimal portfolios with $\gamma_n=10\%$ over 500 simulations based
on $\hSig$ (solid curve) and $\hSig\sam$ (dashed curve) against $p$.
(b) The MSEs of estimated global minimum variances over 500
simulations based on $\hSig$ (solid curve) and $\hSig\sam$ (dashed
curve) against $p$.}}%
\end{center}%
\end{figure}%
\end{singlespace}

\begin{singlespace}
\begin{figure} \centering
\begin{center}%
{\includegraphics[ trim=0.000000in 0.000000in 0.000000in
-0.189256in, height=2in, width=2.6in
]%
{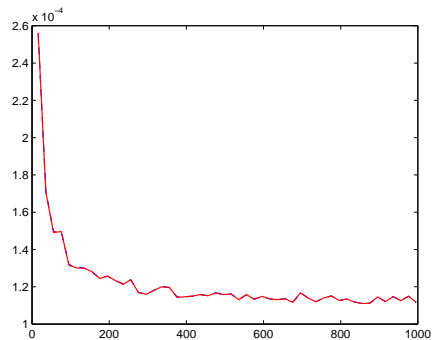}%
}%
\\
\caption{{\scriptsize The MSEs of estimated variances of the equally
weighted portfolio over 500 simulations based on $\hSig$ (solid
curve) and $\hSig\sam$ (dashed curve)
against $p$.}}%
\end{center}%
\end{figure}%
\end{singlespace}

\bigskip
\setcounter{equation}{0} \setcounter{subsect}{0}
\setcounter{subsubsect}{0}

\sect {\bf Concluding remarks.}\quad This paper investigates the
impact of dimensionality on the estimation of covariance matrices.
Two estimators are singled out for studies and comparisons:  the
sample covariance matrix and the factor-model based estimate.  The
inverse of the covariance matrix takes advantage of the factor
structure and hence can be better estimated in the factor approach.
As a result, when the parameters involve the inverse of the
population covariance, substantial gain can be made.  On the other
hand, the covariance matrix itself does not take much advantage of
the factor structure, and hence its estimate can not be improved
much in the factor approach. This is somewhat surprising and is
against the conventional wisdom.

Optimal portfolio allocation and minimum variance portfolio involve
the inverse of the covariance matrix.  Hence, it is advantageous to
employ the factor structure in portfolio allocation.  On the other
hand, intrinsically the risk management does not depend on the
covariance structure and hence there is no advantage to appeal to
the factor model in risk management.

Our conclusion is also verified by an extensive simulation study, in
which the parameters are taken in a neighborhood that is close to
the reality.  The choice of parameters relies on a fit to the famous
Fama-French three-factor model to the portfolios traded in the
market.

Our studies also reveal that the impact of dimensionality on the
estimation of covariance matrices is severe.  This should be taken
into consideration in practical implementations.

\bigskip
\setcounter{equation}{0} \setcounter{subsect}{0}
\setcounter{subsubsect}{0}

\sect {\bf Proofs of theorems.}\quad In this section, we give
rigorous proofs of Theorems 1--7.

\bigskip

\textsc{Proof of Theorem 1}.\quad (1) First, we prove
$\left(pK\right)^{-1}n^{1/2}$-consistency of $\hSig$ under the
Frobenius norm. To facilitate the presentation, we introduce here
some notation used throughout the rest of the paper. Let
$\bC_n\heq\bE\bX\t (\bX\bX\t )^{-1}$,
\[ \bD_n\heq\left\{\left(n-1\right)^{-1}\bX\bX\t
-\left[n(n-1)\right]^{-1}\bX\bone\bone\t \bX\t\right\}-\cov(\bff)\]
and
\[ \bF_n\heq
I_{p}\circ n^{-1}\bE\left(I_n-\bH\right)\bE\t-\Sig_0,
\]
where $\bH\heq\bX\t\left(\bX\bX\t\right)^{-1}\bX$ is the $n\times n$
hat matrix and $\bA_1\circ\bA_2$ stands for the Hadamard product,
i.e. the entrywise product, for any $q\times r$ matrices $\bA_1$ and
$\bA_2$. Then we have
$\hB=\bY\bX\t\left(\bX\bX\t\right)^{-1}=\bB+\bC_n$,
$\hcov(\bff)=\left(n-1\right)^{-1}\bX\bX\t
-\left\{n\left(n-1\right)\right\}^{-1}\bX\bone\bone\t
\bX\t=\cov(\bff)+\bD_n$,
$\hSig_0=\diag\left(n^{-1}\hE\hE\t\right)=\Sig_0+\bF_n$ and
\begin{equation} \label{110}
\hSig=\Sig+\bB\bD_n\bB\t +\left[\bB\hcov(\bff)\bC_n\t
+\bC_n\hcov(\bff)\bB\t\right]+\bC_n\hcov(\bff)\bC_n\t +\bF_n,
\end{equation}
This shows that $\hSig$ is a four-term perturbation of the
population covariance matrix, and this representation is our key
technical tool. By the Cauchy-Schwarz inequality, it follows from
(\ref{110}) that
\begin{align*}
E\|\hSig-\Sig\|^2&\leq 4\ \Big[E\
\tr\left\{\left(\bB\bD_n\bB\t\right)^2\right\}+E\
\tr\left\{\left[\bB\hcov(\bff)\bC_n\t +\bC_n\hcov(\bff)\bB\t\right]^2\right\}\\
&\quad\quad+E\ \tr\left\{\left[\bC_n\hcov(\bff)\bC_n\t
\right]^2\right\}+E\ \tr\left(\bF_n^2\right)\Big].
\end{align*}
We will examine each of the above four terms on the right hand side
separately. For brevity of notation, we suppress the first subscript
$n$ in some situations where the dependence on $n$ is self-evident.

Before going further, let us bound $\|\bB_n\|$. From assumption (B),
we know that $\cov(\bff)\geq\sigma_1I_K$, where for any symmetric
positive semidefinite matrices $\bA_1$ and $\bA_2$, $\bA_1\geq\bA_2$
means $\bA_1-\bA_2$ is positive semidefinite. Thus it follows easily
from (\ref{102}) that
\[ \sigma_1\bB_n\bB_n\t=\bB_n\left(\sigma_1I_K\right)\bB_n\t\leq\bB_n\cov(\bff)\bB_n\t\leq\Sig_n, \]
which along with $b_n=O(p)$ in assumption (B) shows that
$\|\bB_n\|^2=\tr\left(\bB_n\bB_n\t\right)\leq\tr\left(\Sig_n\right)/\sigma_1\leq\frac{b_n}{\sigma_1}=O(p)$,
i.e.
\begin{equation} \label{002}
\|\bB_n\|=O(p^{1/2}).
\end{equation}
Clearly, $\|\bB_n\t\bB_n\|=\|\bB_n\bB_n\t\|$, and by (\ref{018}) in
Lemma 1 and (\ref{002}) we have
\begin{equation} \label{003}
\|\bB_n\t\bB_n\|=\|\bB_n\bB_n\t\|\leq\|\bB_n\|\|\bB_n\t\|=\|\bB_n\|^2=O(p).
\end{equation}
This fact is a key observation that will be used very often, and as
shown above, it is entailed only by assumptions (A) and (B), which
are valid throughout the paper.

Now we consider the first term, say $E\
\tr\{\left(\bB\bD_n\bB\t\right)^2\}$. From $c_n=O(1)$ in assumption
(B), we see that the fourth moments of $\bff$ are bounded across
$n$, thus a routine calculation reveals that
\begin{equation} \label{005}
E\left(\|\bD_n\|^2\right)=O(n^{-1}K^2),
\end{equation}
which is an important fact that will be used very often and also
helps study the inverse $\hcov(\bff)^{-1}$ by keeping in mind that
$K\rightarrow\infty$. By (\ref{019}) in Lemma 1, (\ref{003}), and
(\ref{005}), we have
\begin{equation} \label{006}
E\
\tr\left[\left(\bB\bD_n\bB\t\right)^2\right]\leq\left\|\bB\t\bB\right\|^2E\left(\|\bD_n\|^2\right)
=O(n^{-1}(pK)^2).
\end{equation}

The remaining three terms are taken care of by Lemmas 2 and 3.
Therefore, in view of (\ref{003}), combining (\ref{006}) with
(\ref{012})--(\ref{112}) in Lemmas 2 and 3 gives
\[
E\left\|\hSig-\Sig\right\|^2=O(n^{-1}(pK)^2).
\]
In particular, this implies that
$\left\|\hSig-\Sig\right\|=O_P(n^{-1/2}pK)$, which proves
$\left(pK\right)^{-1}n^{1/2}$-consistency of the covariance matrix
estimator $\hSig$ under Frobenius norm.

\medskip

(2) Then, we show that $\hSig\sam$ is
$\left(pK\right)^{-1}n^{1/2}$-consistent under the Frobenius norm.
By (\ref{102}) and (\ref{106}), we have
\begin{align} \label{147}
\hSig\sam&=\Sig+\bB\bD_n\bB\t+\bG_n+\left(n-1\right)^{-1}\left\{\bB\bX\bE\t+\bE\bX\t\bB\right\}\\
\nonumber&\quad-\left[n\left(n-1\right)\right]^{-1}\left\{\bB\bX\bone\bone\t\bE\t+\bE\bone\bone\t\bX\t\bB\t\right\},
\end{align}
where $\bG_n\heq\left\{\left(n-1\right)^{-1}\bE\bE\t
-\left[n(n-1)\right]^{-1}\bE\bone\bone\t \bE\t\right\}-\Sig_0$. This
shows that $\hSig\sam$ is also a four-term perturbation of the
population covariance matrix. By the Cauchy-Schwarz inequality, it
follows from (\ref{147}) that
\begin{align*}
E\left\|\hSig\sam-\Sig\right\|^2&\leq4\
\Big[E\left\|\bB\bD_n\bB\t\right\|^2+E\left\|\bG_n\right\|^2+
2\left(n-1\right)^{-2}E\left\|\bB\bX\bE\t\right\|^2\\
&\quad\quad+2\left[n\left(n-1\right)\right]^{-2}E\left\|\bB\bX\bone\bone\t\bE\t\right\|^2\Big].
\end{align*}
As in part (1), we will examine each of the above four terms on the
right hand side separately. The first term
$E\left\|\bB\bD_n\bB\t\right\|^2$ has been bounded in (\ref{006}).
Using the same argument as in Lemma 6, we can show that
$E\left\|\bG_n\right\|^2=O(n^{-1}p^2)$. In view of (\ref{003}), it
is shown that
\[ E\left\|\bB\bX\bE\t\right\|^2=O(np^2K) \]
in the proof of Lemma 2. Using the same argument as in Lemma 2 to
bound $E\left\|\bB\bX\bone\bone\t\bH\bE\t\right\|^2$, we can easily
get
\[
E\left\|\bB\bX\bone\bone\t\bE\t\right\|^2=O(n^3p^2K),
\]
which along with (\ref{006}) and the above results yields
\[ E\left\|\hSig\sam-\Sig\right\|^2=O(n^{-1}(pK)^2). \] This proves
$\left(pK\right)^{-1}n^{1/2}$-consistency of $\hSig\sam$ under the
Frobenius norm.

\medskip

(3) Finally, we prove the uniform weak convergence of eigenvalues.
It follows from Corollary 6.3.8 of Horn and Johnson (1985) that
\[ \max_{1\leq k\leq
p}\left|\lambda_k(\hSig_n)-\lambda_k(\Sig_n)\right|\leq\left\{\sum_{k=1}^{p}\left[\lambda_k(\hSig_n)-
\lambda_k(\Sig_n)\right]^2\right\}^{1/2}\leq\left\|\hSig_n-\Sig_n\right\|.
\]
Therefore, the uniform weak convergence of the eigenvalues of the
$\hSig_n$'s follows immediately from the
$\left(pK\right)^{-1}n^{1/2}$-consistency of $\hSig$ under the
Frobenius norm shown in part (1). Similarly, by the
$\left(pK\right)^{-1}n^{1/2}$-consistency of $\hSig\sam$ under the
Frobenius norm shown in part (2), the same conclusion holds for
$\hSig\sam$. $\quad\square$

\bigskip

\textsc{Proof of Theorem 2}.\quad (1) First, we show that $\hSig$ is
$n^{\beta/2}$-consistent under norm $\|\cdot\|_\Sig$. The main idea
of the proof is similar to that of Theorem 1, but the proof is more
tricky and involved here since the norm $\|\cdot\|_\Sig$ involves
the inverse of the covariance matrix $\Sig$. By the Cauchy-Schwarz
inequality, it follows from (\ref{110}) that
\begin{align*}
E\left\|\hSig-\Sig\right\|_\Sig^2&\leq 4\ \Big[E\left\|\bB\bD_n\bB\t\right\|_\Sig^2+
E\left\|\bB\hcov(\bff)\bC_n\t +\bC_n\hcov(\bff)\bB\t\right\|_\Sig^2\\
&\quad\quad+E\left\|\bC_n\hcov(\bff)\bC_n\t\right\|_\Sig^2\}+E\left\|\bF_n\right\|_\Sig^2\Big].
\end{align*}
As in the proof of Theorem 1, we will study each of the above four
terms on the right hand side separately.

Before going further, let us bound
$\left\|\bB\t\Sig^{-1}\bB\right\|$. From (\ref{102}), we know that
$\Sig=\Sig_0+\bB\cov(\bff)\bB\t$, which along with the
Sherman-Morrison-Woodbury formula shows that
\begin{equation} \label{109}
\Sig^{-1}=\Sig_0^{-1}-\Sig_0^{-1}\bB\left[\cov(\bff)^{-1}+\bB\t\Sig_0^{-1}\bB\right]^{-1}\bB\t\Sig_0^{-1}.
\end{equation}
Thus it follows that
\begin{align*}
\bB\t\Sig^{-1}\bB&=\bB\t\Sig_0^{-1}\bB-\bB\t\Sig_0^{-1}\bB\left[\cov(\bff)^{-1}+
\bB\t\Sig_0^{-1}\bB\right]^{-1}\bB\t\Sig_0^{-1}\bB\\
&=\bB\t\Sig_0^{-1}\bB\left[\cov(\bff)^{-1}+\bB\t\Sig_0^{-1}\bB\right]^{-1}\cov(\bff)^{-1}\\
&=\cov(\bff)^{-1}-\cov(\bff)^{-1}\left[\cov(\bff)^{-1}+\bB\t\Sig_0^{-1}\bB\right]^{-1}\cov(\bff)^{-1},
\end{align*}
which implies that
\[ \left\|\bB\t\Sig^{-1}\bB\right\|\leq\left\|\cov(\bff)^{-1}\right\|+
\left\|\cov(\bff)^{-1}\left[\cov(\bff)^{-1}+\bB\t\Sig_0^{-1}\bB\right]^{-1}\cov(\bff)^{-1}\right\|.
\]
Note that $\cov(\bff)^{-1}$ is symmetric positive definite and
$\bB\t\Sig_0^{-1}\bB$ is symmetric positive semidefinite. Thus,
$\cov(\bff)^{-1}+\bB\t\Sig_0^{-1}\bB\geq\cov(\bff)^{-1}$, which in
turn implies that
$\left[\cov(\bff)^{-1}+\bB\t\Sig_0^{-1}\bB\right]^{-1}\leq\cov(\bff)$
and
\[ \cov(\bff)^{-1}\left[\cov(\bff)^{-1}+\bB\t\Sig_0^{-1}\bB\right]^{-1}\cov(\bff)^{-1}\leq
\cov(\bff)^{-1}\cov(\bff)\cov(\bff)^{-1}=\cov(\bff)^{-1}.
\]
In particular, this entails that
\[ \left\|\cov(\bff)^{-1}\left[\cov(\bff)^{-1}+\bB\t\Sig_0^{-1}\bB\right]^{-1}\cov(\bff)^{-1}\right\|\leq
\left\|\cov(\bff)^{-1}\right\|, \] so now the problem of bounding
$\left\|\bB\t\Sig^{-1}\bB\right\|$ reduces to bounding
$\left\|\cov(\bff)^{-1}\right\|$. By assumption (B),
$\lambda_K(\cov(\bff))\geq\sigma_1$ for some constant $\sigma_1>0$.
Thus the largest eigenvalues of $\cov(\bff)^{-1}$ are bounded across
$n$, which easily implies that
$\left\|\cov(\bff)^{-1}\right\|=O(K^{1/2})$. This together with the
above results shows that
\begin{equation} \label{148}
\left\|\bB\t\Sig^{-1}\bB\right\|=O(K^{1/2}).
\end{equation}

Now we are ready to examine the first term, say
$E\left\|\bB\bD_n\bB\t\right\|_\Sig^2$. By (\ref{018}) in Lemma 1,
we have
\[
\left\|\bB\bD_n\bB\t\right\|_\Sig^2=p^{-1}\tr\left[\left(\bD_n\bB\t\Sig^{-1}\bB\right)^2\right]
\leq p^{-1}\left\|\bD_n\right\|^2\left\|\bB\t\Sig^{-1}\bB\right\|^2.
\]
Therefore, it follows from (\ref{005}) and (\ref{148}) that
\begin{equation} \label{111}
E\left\|\bB\bD_n\bB\t\right\|_\Sig^2=O(n^{-1}p^{-1}K^3).
\end{equation}

Then, we consider the second term $E\left\|\bB\hcov(\bff)\bC_n\t
+\bC_n\hcov(\bff)\bB\t\right\|_\Sig^2$. Note that
\begin{align} \label{117}
E\big\|\bB\hcov(\bff)\bC_n\t
&+\bC_n\hcov(\bff)\bB\t\big\|_\Sig^2\leq2\left[E\left\|\bB\hcov(\bff)\bC_n\t\right\|_\Sig^2
+E\left\|\bC_n\hcov(\bff)\bB\t\right\|_\Sig^2\right]\\
\nonumber &=4\ E\left\|\bB\hcov(\bff)\bC_n\t\right\|_\Sig^2\leq8\Big[\left(n-1\right)^{-2}E\left\|\bB\bX\bX\t\bC_n\t\right\|_\Sig^2\\
\nonumber &\quad+n^{-2}\left(n-1\right)^{-2}E\left\|\bB\bX\bone\bone\t\bX\t\bC_n\t\right\|_\Sig^2\Big]\\
\nonumber
&\heq8\left(n-1\right)^{-2}\mathcal{L}_1+8n^{-2}\left(n-1\right)^{-2}\mathcal{L}_2.
\end{align}
Since $E(\bveps|\bff)=\bzero$, conditioning on $\bX$ gives
\begin{align*}
\mathcal{L}_1&=p^{-1}E\ \tr\left[\bX
E\left(\bE\t\Sig^{-1}\bE|\bX\right)\bX\t\bB\t\Sig^{-1}\bB\right]\\
&=p^{-1}E\
\tr\left[\bX\ \tr\left(\Sig^{-1}\Sig_0\right)I_n\ \bX\t\bB\t\Sig^{-1}\bB\right]\\
&\leq
p^{-1}\tr\left(\Sig^{-1}\Sig_0\right)E\left(\|\bX\bX\t\|\right)\left\|\bB\t\Sig^{-1}\bB\right\|.
\end{align*}
In the proof of Lemma 2, it is shown that
$E\left(\|\bX\bX\t\|^2\right)=O(n^2K^2)$, which implies that
\[ E\left(\|\bX\bX\t\|\right)\leq\left[E\left(\|\bX\bX\t\|^2\right)\right]^{1/2}=O(nK). \]
By (\ref{102}) and assumptions (B) and (C), we can easily get
\[ \tr\left(\Sig^{-1}\Sig_0\right)\leq\tr\left(\Sig^{-1}\right)O(1)=O(p), \]
which along with (\ref{148}) and the above results shows that
\[ \mathcal{L}_1=O(nK^{3/2}). \]

Similarly, by conditioning on $\bX$ we have
\begin{align*}
\mathcal{L}_2&=p^{-1}E\ \tr\left[\bX\bone\bone\t\bH
E\left(\bE\t\Sig^{-1}\bE|\bX\right)\bH\bone\bone\t\bX\t\bB\t\Sig^{-1}\bB\right]\\
&=p^{-1}E\ \tr\left[\bX\bone\bone\t\bH \
\tr\left(\Sig^{-1}\Sig_0\right)I_n\
\bH\bone\bone\t\bX\t\bB\t\Sig^{-1}\bB\right].
\end{align*}
Then, applying (\ref{018})--(\ref{020}) in Lemma 1 gives
\begin{align*}
\mathcal{L}_2&\leq
p^{-1}\tr\left(\Sig^{-1}\Sig_0\right)E\left\|\bX\bone\bone\t\bH\bone\bone\t\bX\t\right\|\left\|\bB\t\Sig^{-1}\bB\right\|\\
&\leq
p^{-1}\tr\left(\Sig^{-1}\Sig_0\right)E\left\|\bH\right\|\left\|\bX\t\bX\right\|\left\|\bone\bone\t\bone\bone\t\right\|\left\|\bB\t\Sig^{-1}\bB\right\|\\
&=n^2p^{-1}K^{1/2}\tr\left(\Sig^{-1}\Sig_0\right)E\left\|\bX\t\bX\right\|\left\|\bB\t\Sig^{-1}\bB\right\|,
\end{align*}
which together with the above results shows that
\[ \mathcal{L}_2=O(n^3K^2). \]
Thus, in view of (\ref{117}) we have
\begin{equation} \label{116}
E\left\|\bB\hcov(\bff)\bC_n\t
+\bC_n\hcov(\bff)\bB\t\right\|_\Sig^2=O(n^{-1}K^2).
\end{equation}

The third and fourth terms are examined in Lemmas 4 and 5,
respectively. Since $K\leq p$ by assumption (A), combining
(\ref{111}) and (\ref{116}) with (\ref{138}) and (\ref{139}) in
Lemmas 4 and 5 results in
\[
E\left\|\hSig-\Sig\right\|_\Sig^2=O(n^{-1}K^2)+O(n^{-2}pK).
\]
In particular, when $K=O(n^{\alpha_1})$ and $p=O(n^\alpha)$ for some
$0\leq\alpha_1<1/2$ and $0\leq\alpha<2-\alpha_1$, we have
\[ \left\|\hSig-\Sig\right\|_\Sig=O_P(n^{-\beta/2}) \] with
$\beta=\min\left(1-2\alpha_1,2-\alpha-\alpha_1\right)$, which proves
$n^{\beta/2}$-consistency of covariance matrix estimator $\hSig$
under norm $\|\cdot\|_\Sig$.

\medskip

(2) Then, we prove the $n^{\beta_1/2}$-consistency of $\hSig\sam$
under norm $\|\cdot\|_\Sig$. By the Cauchy-Schwarz inequality, it
follows from (\ref{147}) that
\begin{align*}
E\left\|\hSig\sam-\Sig\right\|_{\Sig}^2&\leq4\
\Big[E\left\|\bB\bD_n\bB\t\right\|_{\Sig}^2+E\left\|\bG_n\right\|_{\Sig}^2+
2\left(n-1\right)^{-2}E\left\|\bB\bX\bE\t\right\|_{\Sig}^2\\
&\quad\quad+2\left[n\left(n-1\right)\right]^{-2}E\left\|\bB\bX\bone\bone\t\bE\t\right\|_{\Sig}^2\Big].
\end{align*}
As in part (1), we will examine each of the above four terms on the
right hand side separately. The first term
$E\left\|\bB\bD_n\bB\t\right\|_{\Sig}^2$ has been bounded in
(\ref{111}), and the second term $E\left\|\bG_n\right\|_{\Sig}^2$ is
considered in Lemma 6. The third term
$E\left\|\bB\bX\bE\t\right\|_{\Sig}^2$ is exactly $\mathcal{L}_1$ in
part (1) above. Using the same argument that was used in part (1) to
prove $\mathcal{L}_2$, we can easily get
\[
E\left\|\bB\bX\bone\bone\t\bE\t\right\|_{\Sig}^2=O(n^3K^{3/2}).
\]
Thus, by (\ref{111}) and (\ref{143}) in Lemma 6 along with the above
results, we have
\[ E\left\|\hSig\sam-\Sig\right\|_\Sig^2=O(n^{-1}p^{-1}K^3)+O(n^{-1}p)+O(n^{-1}K^{3/2}). \]
In particular, when $K=O(n^{\alpha_1})$ and $p=O(n^\alpha)$ for some
$0\leq\alpha<1$ and $0\leq\alpha_1<\left(1+\alpha\right)/3$, we have
\[ \left\|\hSig\sam-\Sig\right\|_\Sig=O_P(n^{-\beta_1/2}) \] with
$\beta_1=1-\max(\alpha,3\alpha_1/2,3\alpha_1-\alpha)$, which shows
$n^{\beta_1/2}$-consistency of $\hSig\sam$ under norm
$\|\cdot\|_\Sig$. $\quad\square$

\bigskip

\textsc{Proof of Theorem 3}.\quad (1) First, we prove the weak
convergence of $\hSig\sam^{-1}$ under the Frobenius norm. Note that
$\hSig\sam$ involves sample covariance matrix estimation of
$\Sig_0$, so the technique in part (2) below does not help. In
general, the only available way is as follows. We define
$\bQ_n=\hSig\sam-\Sig_n$. It is a basic fact in matrix theory that
\begin{equation} \label{040}
\left\|\hSig\sam^{-1}-\Sig_n^{-1}\right\|\leq\left\|\Sig_n^{-1}\right\|\frac{\left\|\Sig_n^{-1}\bQ_n\right\|}
{1-\left\|\Sig_n^{-1}\bQ_n\right\|}\leq
\frac{\left\|\Sig_n^{-1}\right\|^2\left\|\bQ_n\right\|}{1-\left\|\Sig_n^{-1}\right\|\left\|\bQ_n\right\|}
\end{equation}
whenever $\left\|\Sig_n^{-1}\right\|\left\|\bQ_n\right\|<1$. From
Theorem 1, we know that \[ \left\|\bQ_n\right\|=O_P(n^{-1/2}pK).
\] By (\ref{141}), we have
$\left\|\Sig_n^{-1}\right\|=O(p^{1/2})$. Since $pK^{1/2}=o((n/\log
n)^{1/4})$ we see that \[
\left\|\Sig_n^{-1}\right\|\left\|\bQ_n\right\|\toP0\quad
\text{and}\quad\sqrt{np^{-4}K^{-2}/\log n}\
\left\|\Sig_n^{-1}\right\|^2\left\|\bQ_n\right\|\toP0. \] It follows
easily that
\[ \sqrt{np^{-4}K^{-2}/\log n}\ \frac{\left\|\Sig_n^{-1}\right\|^2\left\|\bQ_n\right\|}{1-\left\|\Sig_n^{-1}\right\|\left\|\bQ_n\right\|}\toP0, \]
which along with (\ref{040}) shows that
\[ \sqrt{np^{-4}K^{-2}/\log n}\ \left\|\hSig\sam^{-1}-\Sig_n^{-1}\right\|\toP0\quad\text{as
}n\rightarrow\infty. \]

(2) Then, we show the weak convergence of $\hSig^{-1}$ under the
Frobenius norm. The basic idea is to examine the estimation error
for each term of $\hSig^{-1}$, which has an explicit form thanks to
the factor structure. From (\ref{105}), we know that
$\hSig=\hB\hcov(\bff)\hB\t +\hSig_0$, which along with the
Sherman-Morrison-Woodbury formula shows that
\begin{equation} \label{120}
\hSig^{-1}=\hSig_0^{-1}-\hSig_0^{-1}\hB\left[\hcov(\bff)^{-1}+\hB\t\hSig_0^{-1}\hB\right]^{-1}\hB\t\hSig_0^{-1}.
\end{equation}
Thus by (\ref{109}), we have
\begin{align} \label{121}
\left\|\hSig^{-1}-\Sig^{-1}\right\|&\leq\left\|\hSig_0^{-1}-\Sig_0^{-1}\right\|+\left\|\left(\hSig_0^{-1}-\Sig_0^{-1}\right)
\hB\left[\hcov(\bff)^{-1}+\hB\t\hSig_0^{-1}\hB\right]^{-1}\hB\t\hSig_0^{-1}\right\|\\
\nonumber &\quad+\left\|\Sig_0^{-1}
\hB\left[\hcov(\bff)^{-1}+\hB\t\hSig_0^{-1}\hB\right]^{-1}\hB\t\left(\hSig_0^{-1}-\Sig_0^{-1}\right)\right\|\\
\nonumber &\quad+\left\|\Sig_0^{-1}
\left(\hB-\bB\right)\left[\hcov(\bff)^{-1}+\hB\t\hSig_0^{-1}\hB\right]^{-1}\hB\t\Sig_0^{-1}\right\|\\
\nonumber &\quad+\left\|\Sig_0^{-1}
\bB\left[\hcov(\bff)^{-1}+\hB\t\hSig_0^{-1}\hB\right]^{-1}\left(\hB\t-\bB\t\right)\Sig_0^{-1}\right\|\\
\nonumber &\quad+\left\|\Sig_0^{-1}
\bB\left\{\left[\hcov(\bff)^{-1}+\hB\t\hSig_0^{-1}\hB\right]^{-1}-\left[\cov(\bff)^{-1}+
\bB\t\Sig_0^{-1}\bB\right]^{-1}\right\}\bB\t\Sig_0^{-1}\right\|\\
\nonumber
&\heq\mathcal{K}_1+\mathcal{K}_2+\mathcal{K}_3+\mathcal{K}_4+\mathcal{K}_5+\mathcal{K}_6.
\end{align}
To study $\left\|\hSig^{-1}-\Sig^{-1}\right\|$, we need to examine
each of the above six terms $\mathcal{K}_1,\cdots,\mathcal{K}_6$
separately, so it would be lengthy work to check all the details
here. Therefore, we only sketch the idea of the proof and leave the
details to the reader.

From assumption (C), we know that the diagonal entries of $\Sig_0$
are bounded away from 0. Note that $\hSig_0$ and $\Sig_0$ are both
diagonal, and thus, by the same argument as in Lemma 5, we can
easily show that
\begin{equation} \label{122}
\mathcal{K}_1=\left\|\hSig_0^{-1}-\Sig_0^{-1}\right\|=O_P(n^{-1/2}p^{1/2})+O_P(n^{-1}pK^{1/2})=O_P(n^{-1/2}p^{1/2}),
\end{equation}
since $pK^{1/2}=o((n/\log n)^{1/2})$. Now we consider the second
term $\mathcal{K}_2$. By (\ref{018}) in Lemma 1, we have
\begin{align*}
\mathcal{K}_2&\leq\left\|\left(\hSig_0^{-1}-\Sig_0^{-1}\right)\hSig_0^{1/2}\right\|\left\|\hSig_0^{-1/2}
\hB\left[\hcov(\bff)^{-1}+\hB\t\hSig_0^{-1}\hB\right]^{-1}\hB\t\hSig_0^{-1/2}\right\|\left\|\hSig_0^{-1/2}\right\|\\
&\heq\mathcal{L}_1\mathcal{L}_2\left\|\hSig_0^{-1/2}\right\|,
\end{align*}
and we will examine each of the above two terms $\mathcal{L}_1$ and
$\mathcal{L}_2$, as well as $\left\|\hSig_0^{-1/2}\right\|$. Since
$\hSig_0$ and $\Sig_0$ are diagonal, a similar argument to that
bounding $\mathcal{K}_1$ above applies to show that
\[ \left\|\hSig_0^{-1/2}\right\|=O_P(p^{1/2})\quad\text{and}\quad\mathcal{L}_1=O_P(n^{-1/2}p^{1/2}). \]
Clearly,
$\hSig_0^{-1/2}\hB\left[\hcov(\bff)^{-1}+\hB\t\hSig_0^{-1}\hB\right]^{-1}\hB\t\hSig_0^{-1/2}$
is symmetric positive semidefinite with rank at most $K$ and
$\hSig_0^{1/2}\hSig^{-1}\hSig_0^{1/2}\geq0$. Thus it follows from
(\ref{120}) that
\[ \hSig_0^{-1/2}
\hB\left[\hcov(\bff)^{-1}+\hB\t\hSig_0^{-1}\hB\right]^{-1}\hB\t\hSig_0^{-1/2}=I_p-\hSig_0^{1/2}\hSig^{-1}\hSig_0^{1/2}
\leq I_p,
\]
which implies that
$\hSig_0^{-1/2}\hB\left[\hcov(\bff)^{-1}+\hB\t\hSig_0^{-1}\hB\right]^{-1}\hB\t\hSig_0^{-1/2}$
has at most $K$ positive eigenvalues and all of them are bounded by
one. This shows that $\mathcal{L}_2\leq K^{1/2}$, which along with
the above results gives
\begin{equation} \label{123}
\mathcal{K}_2=O_P(n^{-1/2}pK^{1/2}).
\end{equation}
Similarly, we can also show that
\begin{equation} \label{124}
\mathcal{K}_3=O_P(n^{-1/2}pK^{1/2}).
\end{equation}

Then we consider terms $\mathcal{K}_4$ and $\mathcal{K}_5$. Clearly,
$\hcov(\bff)^{-1}+\hB\t\hSig_0^{-1}\hB\geq\hcov(\bff)^{-1}$, which
in turn entails that
$\left[\hcov(\bff)^{-1}+\hB\t\hSig_0^{-1}\hB\right]^{-1}\leq\hcov(\bff)$
and
\[ \left\|\left[\hcov(\bff)^{-1}+\hB\t\hSig_0^{-1}\hB\right]^{-1}\right\|\leq\left\|\hcov(\bff)\right\|. \]
It is easy to show that $\left\|\hcov(\bff)\right\|=O_P(K)$. Thus we
have
\begin{align} \label{125}
\mathcal{K}_4&\leq\left\|\Sig_0^{-1}
\left(\hB-\bB\right)\right\|\left\|\left[\hcov(\bff)^{-1}+\hB\t\hSig_0^{-1}\hB\right]^{-1}\right\|\left\|\hB\t\Sig_0^{-1}\right\|\\
\nonumber &=O_P(n^{-1}p^{1/2})O_P(K)O_P(p^{1/2})=O_P(n^{-1/2}pK)
\end{align}
and
\begin{align} \label{126}
\mathcal{K}_5&\leq\left\|\Sig_0^{-1}
\bB\right\|\left\|\left[\hcov(\bff)^{-1}+\hB\t\hSig_0^{-1}\hB\right]^{-1}\right\|\left\|\left(\hB\t-\bB\t\right)\Sig_0^{-1}\right\|\\
\nonumber &=O_P(p^{1/2})O_P(K)O_P(n^{-1}p^{1/2}K)=O_P(n^{-1/2}pK).
\end{align}
Finally, by the same argument as in part (1) above, we can show that
\[ \left\|\left[\hcov(\bff)^{-1}+\hB\t\hSig_0^{-1}\hB\right]^{-1}-\left[\cov(\bff)^{-1}
+\bB\t\Sig_0^{-1}\bB\right]^{-1}\right\|=o_P(\left(n/\log
n\right)^{-1/2}K^2).
\] Thus by (\ref{019}) in Lemma 1, we have
\begin{align} \label{127}
\mathcal{K}_6&\leq\left\|\left[\hcov(\bff)^{-1}+\hB\t\hSig_0^{-1}\hB\right]^{-1}-\left[\cov(\bff)^{-1}
+\bB\t\Sig_0^{-1}\bB\right]^{-1}\right\|\left\|\bB\t\Sig_0^{-2}
\bB\right\|\\
\nonumber &=o_P(\left(n/\log
n\right)^{-1/2}K^2)O(p)=o_P(\left(n/\log n\right)^{-1/2}pK^2).
\end{align}
Therefore, it follows from (\ref{121})--(\ref{127}) that
\[ \sqrt{np^{-2}K^{-4}/\log n}\ \left\|\hSig_n^{-1}-\Sig_n^{-1}\right\|\toP0\quad\text{as
}n\rightarrow\infty, \] which completes the proof. $\quad\square$

\bigskip

\textsc{Proof of Theorem 4}.\quad We aim at establishing asymptotic
normality of the $K\times K$ matrix
$\sqrt{n}p^{-2}\bB\t\left(\hSig-\Sig\right)\bB$, and only here are
the $K$ factors $f_1,\cdots,f_K$ assumed fixed across $n$. The basic
idea is to use its four-term decomposition below and to show that
the first term has asymptotic normality by the classical central
limit theorem, while the remaining three terms are all negligible,
say $o_P(1)$, which along with Slutsky's theorem leads to the
desired conclusion. In view of (\ref{110}), we have
\begin{align} \label{024}
\nonumber
\sqrt{n}p^{-2}\bB\t\left(\hSig-\Sig\right)\bB&=\sqrt{n}p^{-2}\bB\t\bB\bD_n\bB\t\bB
+\sqrt{n}p^{-2}\bB\t\left\{\bB\hcov(\bff)\bC_n\t
+\bC_n\hcov(\bff)\bB\t
\right\}\bB\\
\nonumber &\quad+\sqrt{n}p^{-2}\bB\t\bC_n\hcov(\bff)\bC_n\t\bB
+\sqrt{n}p^{-2}\bB\t\bF_n\bB\\
&\heq\mathcal{A}_1+\mathcal{A}_2+\mathcal{A}_3+\mathcal{A}_4.
\end{align}
We will study each of the above four terms
$\mathcal{A}_1,\cdots,\mathcal{A}_4$ separately.

First, we consider the term $\mathcal{A}_1$. Define
\[
\mathcal{H}_n=\frac{n}{n-1}\left(n^{-1}\sum_{i=1}^n\bff_i-E\bff\right)\left(n^{-1}\sum_{i=1}^n\bff_i\t
-E\bff\t\right). \] Then we have
\begin{equation} \label{128}
\hcov(\bff)=\left(n-1\right)^{-1}\sum_{i=1}^n\left(\bff_i-E\bff\right)\left(\bff_i\t
-E\bff\t\right)-\mathcal{H}_n.
\end{equation}
By the classical central limit theorem, we know that
\[ \sqrt{n}\left(n^{-1}\sum_{i=1}^n\bff_i-E\bff\right)\toD \mathcal{N}\left(\bzero,\cov(\bff)\right). \]
It follows from the law of large numbers that
$n^{-1}\sum_{i=1}^n\bff_i-E\bff\toP\bzero$. Thus, by Slutsky's
theorem we have $\sqrt{n}\mathcal{H}_n\toD\bzero$, which in turn
implies that \[ \sqrt{n}\mathcal{H}_n\toP\bzero; \] that is,
$\mathcal{H}_n=o_P(n^{-1/2})$. So in view of (\ref{128}), we have
\begin{equation} \label{129}
\hcov(\bff)=n^{-1}\sum_{i=1}^n\left(\bff_i-E\bff\right)\left(\bff_i\t-E\bff\t
\right)+o_P(n^{-1/2}).
\end{equation}
Therefore, it follows easily from $p^{-1}\bB_n\t\bB_n\rightarrow
\bA$ and (\ref{129}) that
\begin{equation} \label{130}
\mathcal{A}_1=\bA\left\{n^{-1/2}\sum_{i=1}^n\left[\left(\bff_i-E\bff\right)\left(\bff_i\t-E\bff\t
\right)-\cov(\bff)\right]\right\}\bA+o_P(1).
\end{equation}
We define \[
n^{-1/2}\sum_{i=1}^n\left[\left(\bff_i-E\bff\right)\left(\bff_i\t-E\bff\t
\right)-\cov(\bff)\right]\heq U_n=(u_{ij})_{K\times K}. \] By the
classical central limit theorem, we know that [see, e.g. Muirhead
(1982)]
\begin{equation} \label{131}
\vech\left(U_n\right)\toD\mathcal{N}\left(0,H\right),
\end{equation}
where $H$ is determined in an obvious way by
\[ \cov\left(u_{ij},u_{kl}\right)=\kappa^{ijkl}+\kappa^{ik}\kappa^{jl}+
\kappa^{il}\kappa^{jk}, \] with $\kappa^{i_1\cdots i_r}$ the central
moment $E\left[(f_{i_1}-Ef_{i_1})\cdots(f_{i_r}-Ef_{i_r})\right]$ of
$\bff=(f_1,\cdots,f_K)\t$. It follows easily from (\ref{130}) and
(\ref{131}) that
\begin{equation} \label{132}
\vech\left(\mathcal{A}_1\right)\toD\mathcal{N}\left(0,G\right),
\end{equation}
where $G=P_D\left(\bA\otimes \bA\right)DHD\t\left(\bA\otimes
\bA\right)P_D\t$, $D$ is the duplication matrix of order $K$, and
$P_D=(D\t D)^{-1}D\t$.

Then, we examine the second term $\mathcal{A}_2$. From
$p^{-1}\bB_n\t\bB_n\rightarrow \bA$, we know that
\begin{equation} \label{133}
\left\|\bB_n\t\bB_n\right\|=\left\|\bB_n\bB_n\t\right\|=O(p),
\end{equation}
which is in line with (\ref{003}). It follows that
\begin{align} \label{134}
\left\|\mathcal{A}_2\right\|&\leq2\left\|\sqrt{n}p^{-2}\bB\t\bB\hcov(\bff)\bC_n\t\bB\right\|
\leq2n^{1/2}p^{-2}\left\|\bB\t\bB\right\|\left\|\hcov(\bff)\bC_n\t\bB\right\|\\
\nonumber &\leq2n^{1/2}p^{-2}\left\|\bB\t\bB\right\|
\Big\{\left(n-1\right)^{-1}\left\|\bX\bE\t\bB\right\|+n^{-1}\left(n-1\right)^{-1}\left\|\bX\bone\bone\t\bH\bE\t\bB\right\|\Big\}\\
\nonumber
&=O(n^{-1/2}p^{-1})\left\|\bX\bE\t\bB\right\|+O(n^{-3/2}p^{-1})\left\|\bX\bone\bone\t\bH\bE\t\bB\right\|.
\end{align}
Since $E(\bveps|\bff)=\bzero$ and $\Sig_0$ is diagonal, conditioning
on $\bX$ gives
\begin{align*}
E\left\|\bX\bE\t\bB\right\|^2&=E\ \tr\left[\bX
E\left(\bE\t\bB\bB\t\bE|\bX\right)\bX\t\right]=E\ \tr\left[\bX\ \tr\left(\bB\bB\t\Sig_0\right)I_n\ \bX\t\right]\\
&=\tr\left(\bB\bB\t\Sig_0\right)E\left\|\bX\right\|^2=O(p)O(n)=O(np).
\end{align*}
Similarly, by conditioning on $\bX$ we have
\begin{align*}
E\left\|\bX\bone\bone\t\bH\bE\t\bB\right\|^2&=E\
\tr\Big[\bX\bone\bone\t\bH
E\left(\bE\t\bB\bB\t\bE\t|\bX\right)\bH\bone\bone\t\bX\t\Big]\\
&=E\ \tr\left[\bX\bone\bone\t\bH \
\tr\left(\bB\bB\t\Sig_0\right)I_n\ \bH\bone\bone\t\bX\t\right]
\end{align*}
and then applying (\ref{019}) and (\ref{020}) in Lemma 1 yields
\begin{align*}
E\left\|\bX\bone\bone\t\bH\bveps\t\bB\right\|^2
&\leq\tr\left(\bB\bB\t\Sig_0\right)E\Big\{
\left\|\bX\t\bX\right\|\left\|\bone\bone\t\bone\bone\t\right\|\left\|\bH\right\|\Big\}\\
&\leq
O(p)n^2K^{1/2}\left\{E\left(\|\bX\t\bX\|^2\right)\right\}^{1/2}=O(n^3p).
\end{align*}
It follows that $\left\|\bX\bE\t\bB\right\|=O_P(n^{1/2}p^{1/2}) $
and $
\left\|\bX\bone\bone\t\bH\bE\t\bB\right\|=O_P(n^{3/2}p^{1/2})$,
 which together with (\ref{134}) shows that
\begin{equation} \label{135}
\mathcal{A}_2=o_P(1);
\end{equation}
that is, $\mathcal{A}_2$ is a negligible term.

Finally, the third and fourth terms $\mathcal{A}_3$ and
$\mathcal{A}_4$ can also be shown to be negligible by invoking Lemma
3. By (\ref{133}) and (\ref{014}) and (\ref{112}) in Lemma 3, we
have
\begin{align*}
E\left\|\bB\t\bC_n\hcov(\bff)\bC_n\t\bB\right\|^2&\leq
\left\|\bB\bB\t\right\|^2E\left\|\bC_n\hcov(\bff)\bC_n\t\right\|^2\\
&=O(p^2)O(n^{-2}p^2)=O(n^{-2}p^4)
\end{align*}
and
\[
E\left\|\bB\t\bF_n\bB\right\|^2\leq\left\|\bB\bB\t\right\|^2E\left\|\bF_n\right\|^2=O(p^2)O(n^{-1}p)=O(n^{-1}p^3).
\]
It follows that
$\left\|\bB\t\bC_n\hcov(\bff)\bC_n\t\bB\right\|=O_P(n^{-1}p^2)$ and
$\left\|\bB\t\bF_n\bB\right\|=O_P(n^{-1/2}p^{3/2})$, which implies
that
\begin{equation} \label{136}
\mathcal{A}_3=o_P(1)\quad\text{and}\quad\mathcal{A}_4=o_P(1).
\end{equation}
Therefore, in view of (\ref{132}), (\ref{135}), and (\ref{136}),
applying Slutsky's theorem gives
\[ \sqrt{n}\ \vech\left[p^{-2}\bB_n\t\left(\hSig_n-\Sig_n\right)\bB_n\right]\toD\mathcal{N}\left(0,G\right), \]
which proves the asymptotic normality of covariance matrix estimator
$\hSig$. $\quad\square$

\bigskip

\textsc{Proof of Theorem 5}.\quad (1) First, we prove the weak
convergence of the estimated global minimum variance based on
$\hSig$. From Theorem 3, we know that
\[ \sqrt{np^{-2}K^{-4}/\log n}\ \left\|\hSig^{-1}-\Sig^{-1}\right\|\toP 0. \]
Note that
\begin{align*}
\left|\widehat{\varphi}_n-\varphi_n\right|&=\left|\bone\t\left(\hSig^{-1}-\Sig^{-1}\right)\bone\right|=
\left|\tr\left[\left(\hSig^{-1}-\Sig^{-1}\right)\bone\bone\t\right]\right|\\
&\leq\left\|\hSig^{-1}-\Sig^{-1}\right\|\left\|\bone\bone\t\right\|=p\left\|\hSig^{-1}-\Sig^{-1}\right\|.
\end{align*}
Thus we have \[ \sqrt{n\left(pK\right)^{-4}/\log n}\
\left|\widehat{\varphi}_n-\varphi_n\right|\toP0.
\] Since all the $\varphi_n$'s are bounded away from zero, it follows
easily that
\[
\sqrt{n\left(pK\right)^{-4}/\log n}\
\left|\hxi_{ng}\t\hSig_n\hxi_{ng}-\bxi_{ng}\t\Sig_n\bxi_{ng}\right|=\sqrt{n\left(pK\right)^{-4}/\log
n}\ \left|\widehat{\varphi}_n^{-1}-\varphi_n^{-1}\right|\toP0.\]

(2) Then, we prove the conclusion for $\hSig\sam$. From Theorem 3,
we know that
\[ \sqrt{np^{-4}K^{-2}/\log n}\ \left\|\hSig\sam^{-1}-\Sig^{-1}\right\|\toP 0. \]
Therefore, the above argument in part (1) applies to show that
\[
\sqrt{np^{-6}K^{-2}/\log n}\
\left|\hxi_{ng}\t\hSig\sam\hxi_{ng}-\bxi_{ng}\t\Sig_n\bxi_{ng}\right|=\sqrt{np^{-6}K^{-2}/\log
n}\
\left|\widehat{\varphi}_n^{-1}-\varphi_n^{-1}\right|\toP0.\quad\square
\]

\bigskip

\textsc{Proof of Theorem 6}.\quad (1) First, we prove the weak
convergence of the estimated variance of the optimal portfolio based
on $\hSig$. From Theorem 3, we know that
\begin{equation} \label{041}
\sqrt{np^{-2}K^{-4}/\log n}\
\left\|\hSig^{-1}-\Sig^{-1}\right\|\toP0,
\end{equation}
and from part (1) in the proof of Theorem 5, we see that
\begin{equation} \label{044}
\sqrt{n\left(pK\right)^{-4}/\log n}\
\left|\widehat{\varphi}_n-\varphi_n\right|\toP0.
\end{equation}
Now we show the same rate for
$\left|\widehat{\psi}_n-\psi_n\right|$, say
\begin{equation} \label{001}
\sqrt{n\left(pK\right)^{-4}/\log n}\
\left|\widehat{\psi}_n-\psi_n\right|\toP0.
\end{equation}
By $b_n=O(p)$ in assumption (B), a routine calculation yields
$\left\|\bmu_n\right\|=O(p^{1/2})$ and
$E\left\|\hmu_n-\bmu_n\right\|^2=O(n^{-1}p)$, and thus
\[ \left\|\hmu_n-\bmu_n\right\|=O_P(n^{-1/2}p^{1/2}). \]
It follows that
\begin{align*}
\left|\widehat{\psi}_n-\psi_n\right|&\leq\left|\bone\t\left(\hSig^{-1}-\Sig^{-1}\right)\hmu\right|+\left|\bone\t\Sig^{-1}
\left(\hmu-\bmu\right)\right|
\leq\left\|\bone\t\right\|\left\|\hSig^{-1}-\Sig^{-1}\right\|\\
&\quad\cdot\left(\left\|\bmu\right\|+\left\|\hmu-\bmu\right\|\right)+\left\|\bone\t\right\|
\left\|\Sig^{-1}\right\|\left\|\hmu-\bmu\right\|.
\end{align*}
Then we have
\begin{align*}
\left|\widehat{\psi}_n-\psi_n\right|&\leq p^{1/2}\left\|\hSig^{-1}-\Sig^{-1}\right\|\left[O(p^{1/2})+O_P(n^{-1/2}p^{1/2})\right]\\
&\quad+p^{1/2}O(p^{1/2})O_P(n^{-1/2}p^{1/2})\\
&=\left\|\hSig^{-1}-\Sig^{-1}\right\|O(p)+O_P(n^{-1/2}p^{3/2})=\left\|\hSig^{-1}-\Sig^{-1}\right\|O(p),
\end{align*}
which together with (\ref{041}) proves (\ref{044}). Similarly, we
can also show that
\begin{equation} \label{046}
\sqrt{n\left(pK\right)^{-4}/\log n}\
\left|\widehat{\phi}_n-\phi_n\right|\toP0.
\end{equation}
Since $\varphi_n\phi_n-\psi_n^2$ are bounded away from zero and
$\varphi_n/(\varphi_n\phi_n-\psi_n^2)$,
$\psi_n/(\varphi_n\phi_n-\psi_n^2)$,
$\phi_n/(\varphi_n\phi_n-\psi_n^2)$, $\gamma_n$ are bounded, the
conclusion follows from (\ref{038}) and (\ref{044})--(\ref{046}).

\medskip

(2) Now we prove the conclusion for $\hSig\sam$. From Theorem 3, we
know that
\[ \sqrt{np^{-4}K^{-2}/\log n}\
\left\|\hSig\sam^{-1}-\Sig^{-1}\right\|\toP0, \] and from part (2)
in the proof of Theorem 5, we see that
\[ \sqrt{np^{-6}K^{-2}/\log n}\
\left|\widehat{\varphi}_n-\varphi_n\right|\toP0. \] Since $b_n=O(p)$
by assumption (B), a routine calculation shows that \[
\left\|\hmu\sam-\bmu_n\right\|=O_P(n^{-1/2}p^{1/2}), \] where
$\hmu\sam$ is the sample mean of $\bmu_n$. Therefore, the argument
in part (1) above applies to show that
\[
\sqrt{np^{-6}K^{-2}/\log n}\
\left|\hxi_n\t\hSig\sam\hxi_n-\bxi_n\t\Sig_n\bxi_n\right|\toP0\quad\text{as
}n\rightarrow\infty.\quad\square \]

\bigskip

\textsc{Proof of Theorem 7}.\quad Since $\bxi_n=O(1)\bone$, the
conclusion follows easily from consistency results of $\hSig$ and
$\hSig\sam$ under the Frobenius norm in Theorem 1. In particular,
when the portfolios $\bxi_n=(\xi_1,\cdots,\xi_{p})\t$ have no short
positions, we have
\[
\left\|\bxi_n\right\|=\sqrt{\xi_1^2+\cdots+\xi_{p}^2}\leq\sqrt{\xi_1+\cdots+\xi_{p}}=1.
\] It therefore follows easily that
\[ \sqrt{n\left(pK\right)^{-2}/\log n}\ \left|\bxi_n\t\hSig_n\bxi_n-\bxi_n\t\Sig_n\bxi_n\right|\toP0\quad\text{as
}n\rightarrow\infty \] and
\[ \sqrt{n\left(pK\right)^{-2}/\log n}\ \left|\bxi_n\t\hSig\sam\bxi_n-\bxi_n\t\Sig_n\bxi_n\right|\toP0\quad\text{as
}n\rightarrow\infty.\quad\square \bigskip \]



\setcounter{equation}{0}
\renewcommand{\theequation}{A.\arabic{equation}}

\begin{center}
APPENDIX
\end{center}

Throughout the paper, we denote by $\bH$ the $n\times n$ hat matrix
$\bX\t\left(\bX\bX\t\right)^{-1}\bX$, which is symmetric and
positive semidefinite with probability one by assumption (A).

\bigskip

\begin{lem} (Basic facts). \em
\begin{itemize}
\item[(i)] For any $q\times r$ matrix $\bA_1$ and $r\times q$ matrix $\bA_2$,
we have
\begin{equation} \label{018}
\left|\tr\left(\bA_1\bA_2\right)\right|\leq\left\|\bA_1\right\|\left\|\bA_2\right\|\
and\
\left\|\bA_1\bA_2\right\|\leq\left\|\bA_1\right\|\left\|\bA_2\right\|.
\end{equation}
In particular, for any $q\times r$ matrix $\bA_1$ and $r\times r$
symmetric matrix $\bA_2$, we have
\begin{equation} \label{019}
\left|\tr\left(\bA_1\bA_2\bA_1\t\right)\right|\leq\left\|\bA_1\t\bA_1\right\|\left\|\bA_2\right\|\
and\
\left\|\bA_1\bA_2\bA_1\t\right\|\leq\left\|\bA_1\t\bA_1\right\|\left\|\bA_2\right\|.
\end{equation}

\item[(ii)] With probability one, the hat matrix $\bH$ is idempotent with
\begin{equation} \label{020}
\tr\left(\bH^2\right)=\tr\left(\bH\right)=K,
\end{equation}
and it satisfies
\begin{equation} \label{021}
0\leq\tr\left(\bH\bone\bone\t\bH\right)\leq
K^{1/2}n\quad\text{and}\quad0\leq\tr\left[\left(\bH\bone\bone\t\bH\right)^2\right]\leq
Kn^2.
\end{equation}
\end{itemize}
\end{lem}

\begin{proof}
One can refer to Horn and Johnson (1990) for standard proofs of
(\ref{018}) and (\ref{019}). The fact that the hat matrix $\bH$ is
idempotent with (\ref{020}) is known in multivariate statistical
analysis. Clearly, $\tr\left(\bH\bone\bone\t\bH\right)=\bone\t
\bH\bone\geq0$. Thus by (\ref{018}) and (\ref{020}), we have
\[
\tr\left(\bH\bone\bone\t\bH\right)=\tr\left(\bH\bone\bone\t\right)
\leq\left\|\bH\right\|\left\|\bone\bone\t\right\|=K^{1/2}n
\]
and
\[
\tr\left[\left(\bH\bone\bone\t\bH\right)^2\right]=
\tr\left[\left(\bH\bone\bone\t\right)^2\right]\leq\left\|\bH\bone\bone\t\right\|^2
\leq\left\|\bH\right\|^2 \left\|\bone\bone\t\right\|^2=Kn^2.
\]
This completes the proof. \quad$\square$
\end{proof}

\bigskip

The main trick in the proofs of the technical lemmas below is
conditioning on $\bX$ and resorting to the basic facts from Lemma 1.

\bigskip

\begin{lem}.\quad\em Under conditions (A) and (B), we have
\begin{equation} \label{012}
E\ \tr\left\{\left[\bB\hcov(\bff)\bC_n\t +\bC_n\hcov(\bff)\bB\t
\right]^2\right\}\leq\left\|\bB\t\bB\right\|O(n^{-1}pK^{3/2}).
\end{equation}
\end{lem}
\begin{proof}
It follows from (\ref{018}) that
\begin{align*}
E\ &\tr\left\{\left[\bB\hcov(\bff)\bC_n\t +\bC_n\hcov(\bff)\bB\t
\right]^2\right\}\leq 2\left(n-1\right)^{-2}E\
\tr\left\{\left[\bB\bX\bE\t+\bE\bX\t\bB\t\right]^2\right\}\\
&\quad+2n^{-2}\left(n-1\right)^{-2}E\
\tr\left\{\left[\bB\bX\bone\bone\t\bH\bE\t+\bE\bH\bone\bone\t\bX\t\bB\t\right]^2\right\}\\
&\heq2\left(n-1\right)^{-2}\mathcal{A}_1+2n^{-2}\left(n-1\right)^{-2}\mathcal{A}_2.
\end{align*}
We will consider the above two terms $\mathcal{A}_1$ and
$\mathcal{A}_2$ separately. By $c_n=O(1)$ in assumption (B), we can
easily get $\left\|E\left(\bff\bff\t\right)\right\|=O(K)$ and
$E\left(\left\|\bff\right\|^4\right)=O(K^2)$.

Since $E(\bE|\bff)=\bzero$, by (\ref{018}) and (\ref{019})
conditioning on $\bX$ results in
\begin{align*}
E\left\|\bB\bX\bE\t\right\|^2&=E\ \tr\left[\bB\bX
E\left(\bE\t\bE|\bX\right)\bX\t\bB\t\right]= E\ \tr\left[\bB\bX\
\tr\left(\Sig_0\right)I_n\
\bX\t\bB\t\right]\\
&=n\ \tr\left(\Sig_0\right)E\ \tr\left[\bB\bff\bff\t\bB\t\right]=n\
\tr\left(\Sig_0\right)\
\tr\left[\bB E\left(\bff\bff\t\right)\bB\t\right]\\
&\leq n\
\tr\left(\Sig_0\right)\left\|\bB\t\bB\right\|\left\|E\left(\bff\bff\t\right)\right\|
=\left\|\bB\t\bB\right\|\tr\left(\Sig_0\right)O(nK).
\end{align*}
Similarly, by conditioning on $\bX$ we have
\begin{align*}
E\left\|\bB\bX\bone\bone\t\bH\bE\t\right\|^2&=E\
\tr\Big[\bB\bX\bone\bone\t\bH
E\left(\bE\t\bE|\bX\right)\bH\bone\bone\t\bX\t\bB\t\Big]\\
&= E\ \tr\left[\bB\bX\bone\bone\t\bH\ \tr\left(\Sig_0\right)I_n\
\bH\bone\bone\t\bX\t\bB\t\right],
\end{align*}
and then applying (\ref{018}) and (\ref{019}) in Lemma 1 gives
\begin{align*}
E\left\|\bB\bX\bone\bone\t\bH\bE\t\right\|^2&\leq
\tr\left(\Sig_0\right)\ E\Big\{
\left\|\bB\t\bB\right\|\left\|\bX\t\bX\right\|\left\|\bone\bone\t\bone\bone\t\right\|
\left\|\bH\right\|\Big\}\\
&\leq K^{1/2}
n^2\left\|\bB\t\bB\right\|\tr\left(\Sig_0\right)\left\{E\left(\left\|\bX\t\bX\right\|^2\right)\right\}^{1/2}.
\end{align*}
Note that $ E\left(\left\|\bX\t\bX\right\|^2\right)=
nE\left(\left\|\bff\right\|^4\right)+n\left(n-1\right)\left\|E\left(\bff\bff\t\right)\right\|^2
=O(n^2K^2)$. Thus,
\[
E\left\|\bB\bX\bone\bone\t\bH\bE\t\right\|^2\leq\left\|\bB\t\bB\right\|\tr\left(\Sig_0\right)O(n^3K^{3/2}).
\] Therefore, by (\ref{018}) we have
\[
\mathcal{A}_1\leq4E\left\|\bB\bX\bE\t\right\|^2
\leq\left\|\bB\t\bB\right\|\tr\left(\Sig_0\right)O(nK) \] and
\[
\mathcal{A}_2\leq
4E\left\|\bB\bX\bone\bone\t\bH\bE\t\right\|^2\leq\left\|\bB\t\bB\right\|\tr\left(\Sig_0\right)O(n^3K^{3/2}),
\]
which together yield (\ref{012}) since clearly
$\tr\left(\Sig_0\right)=O(p)$. $\quad\square$
\end{proof}

\bigskip

\begin{lem}.\quad\em Under conditions (A) and (B), we have
\begin{equation} \label{014}
E\ \tr\left\{\left[\bC_n\hcov(\bff)\bC_n\t
\right]^2\right\}=O(n^{-2}p^2K)
\end{equation}
and
\begin{equation} \label{112}
E\ \tr\left(\bF_n^2\right)=O(n^{-1}pK)+O(n^{-2}p^2K).
\end{equation}
\end{lem}
\begin{proof} The proofs of (\ref{014}) and (\ref{112}) are similar
to those in Lemmas 4 and 5 below, respectively. For brevity, we omit
them here. $\quad\square$
\end{proof}

\bigskip

\begin{lem}.\quad\em Under conditions (A)--(C), we have
\begin{equation} \label{138}
E\left\|\bC_n\hcov(\bff)\bC_n\t\right\|_\Sig^2=O(n^{-2}pK).
\end{equation}
\end{lem}
\begin{proof}
Note that
\begin{align*}
E\left\|\bC_n\hcov(\bff)\bC_n\t\right\|_\Sig^2&\leq2\left(n-1\right)^{-2}E\left\|
\bE\bH\bE\t\right\|_\Sig^2+2n^{-2}\left(n-1\right)^{-2}E\left\|
\bE\bH\bone\bone\t\bH\bE\t\right\|_\Sig^2\\
&\heq2\left(n-1\right)^{-2}\mathcal{K}_1+2n^{-2}\left(n-1\right)^{-2}\mathcal{K}_2.
\end{align*}
We will consider the above two terms $\mathcal{K}_1$ and
$\mathcal{K}_2$ separately. First, we study the term
$\mathcal{K}_1$, which can further be decomposed into four terms.
Since $E\left(\bveps|\bff\right)=\bzero$, by conditioning on $\bX$
we have
\begin{align*}
\mathcal{K}_1&=p^{-1}E\
\tr\left[E\left(\sum_{i,j=1}^n\bH_{ij}\bveps_i\bveps_j\t\Sig^{-1}
\sum_{k,l=1}^n\bH_{kl}\bveps_k\bveps_l\t\Sig^{-1} |\bX\right)\right]\\
&=p^{-1}\mathcal{L}_1+p^{-1}\mathcal{L}_2+p^{-1}\mathcal{L}_3+p^{-1}\mathcal{L}_4,
\end{align*}
where
\begin{align*}
\mathcal{L}_1&=E\
\tr\left\{\sum_{i=1}^n\left(\bH_{ii}\right)^2E\left[\left(
\bveps_i\bveps_i\t\Sig^{-1} \right)^2\right]\right\},\
\mathcal{L}_2=E\ \tr\left\{\sum_{i\neq
j}\bH_{ii}\bH_{jj}E\left(\bveps_i\bveps_i\t\Sig^{-1}
\bveps_j\bveps_j\t\Sig^{-1}\right)\right\},\\
\mathcal{L}_3&=E\ \tr\left\{\sum_{i\neq j}\left(\bH_{ij}\right)^2
E\left[\left(\bveps_i\bveps_j\t\Sig^{-1} \right)^2\right]\right\},\
\mathcal{L}_4=E\ \tr\left[\sum_{i\neq j}\bH_{ij}\bH_{ji}
E\left(\bveps_i\bveps_j\t\Sig^{-1}
\bveps_j\bveps_i\t\Sig^{-1}\right)\right],
\end{align*}
and $\bH_{ij}$ is the $(i,j)$-entry of the $n\times n$ hat matrix
$\bH$. Then we consider each of these four terms separately. By
(\ref{102}) and assumptions (C) and (B), it is easy to see that
\begin{equation} \label{141}
\tr\left(\Sig^{-1}\right)=O(p),\
\left\|\Sig^{-1}\right\|=O(p^{1/2}),\text{ and
}\tr\left[\left(\Sig_0\Sig^{-1}\right)^2\right]=O(p).
\end{equation}
It follows from (\ref{020}) and (\ref{141}) that
\begin{align*}
\mathcal{L}_1&\leq K\ E\left\{\tr\left[\left(\bveps\bveps\t
\Sig^{-1}\right)^2\right]\right\}=K\
E\left[\sum_{i,j=1}^{p}\left(\Sig^{-1}\right)_{ij}\veps_i\veps_j\sum_{k,l=1}^{p}\left(\Sig^{-1}\right)_{kl}\veps_k\veps_l\right]\\
&=K\sum_{i=1}^{p}
\left(\Sig^{-1}\right)_{ii}^2E\left(\veps_i^4\right)+K\sum_{i\neq
j}\left(\Sig^{-1}\right)_{ii}\left(\Sig^{-1}\right)_{jj}E\left(\veps_i^2\right)E\left(\veps_j^2\right)\\
&\quad\quad+2K\sum_{i\neq
j}E\left(\veps_i^2\right)\left(\Sig^{-1}\right)_{ij}E\left(\veps_j^2\right)\left(\Sig^{-1}\right)_{ji}\\
&\leq\left[\tr\left(\Sig^{-1}\right)\right]^2O(K)+\tr\left(\Sig_0\Sig^{-1}\Sig_0\Sig^{-1}\right)O(K)=O(p^2K)
\end{align*}
and
\begin{align*}
\mathcal{L}_2&=E\left\{\sum_{i\neq j}\bH_{ii}\bH_{jj}\ \tr\left[
E\left(\bveps\bveps\t\right)\Sig^{-1}
E\left(\bveps\bveps\t\right)\Sig^{-1}\right]\right\}\\
&\leq\tr\left[\left(\Sig_0\Sig^{-1}\right)^2\right]E\left\{\left[\tr\left(\bH\right)\right]^2\right\}
=O(pK^2).
\end{align*}
Similarly, we have
\begin{align*}
\mathcal{L}_3&\leq K\ \tr\left\{E\left[\left(\bveps\bet\t\Sig^{-1}
\right)^2\right]\right\}=K\ E\left[\sum_{i,j=1}^{p}(\Sig^{-1})_{ij}\
\eta_i\veps_j\sum_{k,l=1}^{p}(\Sig^{-1})_{kl}\
\eta_k\veps_l\right]\\
&=K\sum_{i=1}^{p}
\left(\Sig^{-1}\right)_{ii}^2E\left(\eta_i^2\right)E\left(\veps_i^2\right)+2K\sum_{i\neq
j}E\left(\eta_i^2\right)\left(\Sig^{-1}\right)_{ij}E\left(\veps_j^2\right)\left(\Sig^{-1}\right)_{ji}\\
&\leq\left\|\Sig^{-1}\right\|^2O(K)+\tr\left(\Sig_0\Sig^{-1}\Sig_0\Sig^{-1}\right)O(K)=O(pK)
\end{align*}
and
\begin{align*}
\mathcal{L}_4&\leq K\ \tr\left[E\left(\bveps\bet\t\Sig^{-1}
\bet\bveps\t\Sig^{-1}\right)\right]=K\left[E\left(\bveps\t\Sig^{-1}
\bveps\right)\right]^2=K\left[\sum_{i=1}^{p}(\Sig^{-1})_{ii}E(\veps_i^2)\right]^2\\
&=K\left[\tr\left(\Sig^{-1}\right)O(1)\right]^2=O(p^2K),
\end{align*}
where $\bet=(\eta_1,\cdots,\eta_{p})\t$ is an independent copy of
$\bveps=(\veps_1,\cdots,\veps_{p})\t$. Since $K\leq p$ by assumption
(A), combining $\mathcal{L}_1$, $\mathcal{L}_2$, $\mathcal{L}_3$,
and $\mathcal{L}_4$ together gives
\begin{equation} \label{142}
\mathcal{K}_1=E\left\| \bE\bH\bE\t\right\|_\Sig^2=O(pK).
\end{equation}

Now we consider the second term $\mathcal{K}_2$. By (\ref{021}), the
same calculation as above applies to show that
\[ \mathcal{K}_2=E\left\|
\bE\bH\bone\bone\t\bH\bE\t\right\|_\Sig^2=O(n^2pK). \] Therefore,
combining the above results together yields (\ref{138}).
$\quad\square$
\end{proof}

\bigskip

\begin{lem}.\quad\em Under conditions (A)--(C), we have
\begin{equation} \label{139}
E\|\bF_n\|_{\Sig}^2=O(n^{-1})+O(n^{-2}pK).
\end{equation}
\end{lem}
\begin{proof}
Note that
\[
E\|\bF_n\|_{\Sig}^2\leq2E\left\|I_{p}\circ
n^{-1}\bE\bE\t-\Sig_0\right\|_{\Sig}^2
+2n^{-2}E\left\|I_{p}\circ\bE\bH\bE\t\right\|_{\Sig}^2.
\]
Since $E\left(\bveps\right)=\bzero$ and
$\cov\left(\bveps|\bff\right)=\Sig_0$, we have
\begin{align*}
E\big\|I_{p}\circ&
n^{-1}\bE\bE\t-\Sig_0\big\|_{\Sig}^2=p^{-1}E\left\|n^{-1}\Sig^{-1/2}\left(I_{p}\circ
\bE\bE\t\right)\Sig^{-1/2}-\Sig^{-1/2}\Sig_0\Sig^{-1/2}\right\|_{\Sig}^2\\
&=p^{-1}n^{-1}\left[E\left\|\Sig^{-1/2}\diag\left(\veps_1^2,\cdots,\veps_{p}^2\right)\Sig^{-1/2}\right\|^2-
\left\|\Sig^{-1/2}\Sig_0\Sig^{-1/2}\right\|^2\right]\\
&\leq p^{-1}n^{-1}E\
\tr\left\{\left[\Sig^{-1/2}\diag\left(\veps_1^2,\cdots,\veps_{p}^2\right)\Sig^{-1/2}\right]^2\right\}
\heq p^{-1}n^{-1}\mathcal{L}. \end{align*} It follows from
(\ref{141}) that
\begin{align*}
\mathcal{L}&=\sum_{i,j=1}^{p}E\left[\veps_i^2\left(\Sig^{-1}\right)_{ij}\veps_j^2\left(\Sig^{-1}\right)_{ji}\right]
=\sum_{i=1}^{p}\left(\Sig^{-1}\right)_{ii}^2E\left(\veps^4\right)+
\sum_{i\neq
j}\left(\Sig^{-1}\right)_{ij}^2\left[E\left(\veps^2\right)\right]^2\\
&=\left\|\Sig^{-1}\right\|^2O(1)=O(p),
\end{align*}
which shows that $E\big\|I_{p}\circ
n^{-1}\bE\bE\t-\Sig_0\big\|_{\Sig}^2=O(n^{-1})$. The argument
proving (\ref{142}) in Lemma 4 applies to show that
\[ E\left\|I_{p}\circ\bE\bH\bE\t\right\|_{\Sig}^2=O(pK). \]
Hence, combining the above results together gives (\ref{139}).
$\quad\square$
\end{proof}

\bigskip

\begin{lem}.\quad\em Under conditions (A)--(C), we have
\begin{equation} \label{143}
E\|\bG_n\|_{\Sig}^2=O(n^{-1}p).
\end{equation}
\end{lem}
\begin{proof}
Recall that $\bG_n\heq\left\{\left(n-1\right)^{-1}\bE\bE\t
-\left[n(n-1)\right]^{-1}\bE\bone\bone\t \bE\t\right\}-\Sig_0$, as
defined in part (2) of the proof of Theorem 1. Note that
\begin{align*}
E\|\bG_n\|_{\Sig}^2&\leq3E\left\|n^{-1}\bE\bE\t-\Sig_0\right\|_{\Sig}^2+
3n^{-2}\left(n-1\right)^{-2}E\left\|\bE\bE\t\right\|_{\Sig}^2\\
&\quad+3n^{-2}\left(n-1\right)^{-2}E\left\|\bE\bone\bone\t
\bE\t\right\|_{\Sig}^2.
\end{align*}
From the proofs of $\mathcal{L}_1$ and $\mathcal{L}_4$ in Lemma 4,
we know that
\[E\left\{\tr\left[\left(\bveps\bveps\t
\Sig^{-1}\right)^2\right]\right\}=O(p^2)\quad\text{and}\quad
E\left[\tr\left(\bveps\bet\t\Sig^{-1}
\bet\bveps\t\Sig^{-1}\right)\right]=O(p^2), \] where
$\bet=(\eta_1,\cdots,\eta_{p})\t$ is an independent copy of
$\bveps=(\veps_1,\cdots,\veps_{p})\t$. Thus, we have
\begin{align*}
E&\left\|n^{-1}\bE\bE\t-\Sig_0\right\|_{\Sig}^2=p^{-1}E\left\|n^{-1}\Sig^{-1/2}\bE\bE\t\Sig^{-1/2}
-\Sig^{-1/2}\Sig_0\Sig^{-1/2}\right\|^2\\
&\leq
p^{-1}n^{-1}E\left\|\Sig^{-1/2}\bveps\bveps\t\Sig^{-1/2}\right\|^2=
p^{-1}n^{-1}E\left\{\tr\left[\left(\bveps\bveps\t
\Sig^{-1}\right)^2\right]\right\}=O(n^{-1}p).
\end{align*}

Similarly, it follows that
\begin{align*}
E\left\|\bE\bE\t\right\|_{\Sig}^2&=E\left\|\left(\bveps_1,\cdots,\bveps_n\right)
\left(\bveps_1,\cdots,\bveps_n\right)\t\right\|_{\Sig}^2\leq
np^{-1}E\left\|\Sig^{-1/2}\bveps\bveps\t\Sig^{-1/2}\right\|^2\\
&=np^{-1}E\left\{\tr\left[\left(\bveps\bveps\t
\Sig^{-1}\right)^2\right]\right\}=O(np)
\end{align*}
and
\begin{align*}
E&\left\|\bE\bone\bone\t \bE\t\right\|_{\Sig}^2 \leq
np^{-1}E\left\|\Sig^{-1/2}\bveps\bveps\t\Sig^{-1/2}\right\|^2+
n\left(n-1\right)p^{-1}E\left\|\Sig^{-1/2}\bveps\bet\t\Sig^{-1/2}\right\|^2\\
&\leq np^{-1}E\left\{\tr\left[\left(\bveps\bveps\t
\Sig^{-1}\right)^2\right]\right\}+
n\left(n-1\right)p^{-1}E\left[\tr\left(\bveps\bet\t\Sig^{-1}
\bet\bveps\t\Sig^{-1}\right)\right]=O(n^2p).
\end{align*}
Therefore, combining the above results together proves (\ref{143}).
$\quad\square$
\end{proof}

\bigskip


\begin{center}
REFERENCES
\end{center}

\begin{singlespace}
\begin{itemize}
\item[] \textsc{Aguilar, O.} and \textsc{West, M.} (2000). Bayesian dynamic factor models and
portfolio allocation. {\em Journal of Business and Economic
Statistics} {\bf 18} 338--357.

\item[] \textsc{Bai, J.} (2003). Inferential theory for factor models of large
dimensions. {\em Econometrica} {\bf 71} 135--171.

\item[] \textsc{Boik, R. J.} (2002). Spectral models for covariance matrices.
{\em Biometrika} {\bf 89} 159--182.

\item[] \textsc{Browne, M. W.} (1987). Robustness of statistical inference in factor analysis and related
models. {\em Biometrika} {\bf 74} 375--384.

\item[] \textsc{Browne, M. W.} and \textsc{Shapiro, A.} (1987). Adjustments for kurtosis in factor analysis with elliptically
distributed errors. {\em J. Roy. Statist. Soc. Ser. B} {\bf 49}
346--352.

\item[] \textsc{Campbell, J. Y.}, \textsc{Lo, A. W.} and \textsc{MacKinlay, A. G.} (1997).
{\em The Econometrics of Financial Markets}. Princeton University
Press, New Jersey.

\item[] \textsc{Chamberlain, G.} (1983). Funds, factors and diversification in Arbitrage
Pricing Theory. {\em Econometrica} {\bf 51} 1305--1323.

\item[] \textsc{Chamberlain, G.} and \textsc{Rothschild, M.} (1983). Arbitrage, factor
structure, and mean-variance analysis on large asset markets. {\em
Econometrica} {\bf 51} 1281--1304.

\item[] \textsc{Chiu, T.Y.M.}, \textsc{Leonard, T.} and \textsc{Tsui, K.W.} (1996). The
matrix-logarithm covariance model. {\em J. Amer. Statist. Assoc.}
{\bf 91} 198--210.

\item[] \textsc{Cochrane, J. H.} (2001). {\em Asset Pricing}. Princeton University
Press, New Jersey.

\item[] \textsc{Dempster, A. P.} (1972). Covariance selection.
{\em Biometrics} {\bf 28} 157--175.

\item[] \textsc{Diebold, F. X.} and \textsc{Nerlove, M.} (1989). The dynamics
of exchange rate volatility: a multivariate latent-factor ARCH
model. {\em Journal of Applied Econometrics} {\bf 4} 1--22.

\item[] \textsc{Diggle, P. J.} and \textsc{Verbyla, A. P.} (1998). Nonparametric estimation of
covariance structure in longitudinal data. {\em Biometrics} {\bf 54}
401--415.

\item[] \textsc{Donoho, D. L.} (2000). High-dimensional data analysis: The curses and blessings of
dimensionality. Aide-Memoire of a Lecture at AMS Conference on Math
Challenges of the 21st Century.

\item[] \textsc{Eaton, M. L.} and \textsc{Tyler, D. E.} (1991). On Wielandt's inequality and
its application to the asymptotic distribution of the eigenvalues of
a random symmetric matrix. {\em Ann. Statist.} {\bf 19} 260--271.

\item[] \textsc{Eaton, M. L.} and \textsc{Tyler, D. E.} (1994). The asymptotic distribution of
singular values with applications to canonical correlations and
corresponding analysis. {\em J. Multiv. Anal.} {\bf 50} 238--264.

\item[] \textsc{Engle, R. F.} and \textsc{Watson, M. W.} (1981). A one-factor multivariate time series model of metropolitan
wage rates. {\em J. Amer. Statist. Assoc.} {\bf 76} 774--781.

\item[] \textsc{Fama, E.} and \textsc{French, K.} (1992). The cross-section of expected stock
returns. {\em Jour. Fin.} {\bf 47} 427--465.

\item[] \textsc{Fama, E.} and \textsc{French, K.} (1993). Common risk factors in the
returns on stocks and bonds. {\em Jour. Fin. Econ.} {\bf 33} 3--56.

\item[] \textsc{Fan, J.} (2005). A selective overview of nonparametric methods in
financial econometrics (with discussion). {\em Statist. Science}, to
appear.

\item[] \textsc{Fan, J.} and \textsc{Li, R.} (2006). Statistical
challenges with high dimensionality: feature selection in knowledge
discovery. {\em Proceedings of the Madrid International Congress of
Mathematicians 2006}, to apepar.

\item[] \textsc{Fan, J.} and \textsc{Peng, H.} (2004). Nonconcave
penalized likelihood with a diverging number of parameters. {\em
Ann. Statist.} {\bf 32} 928--961.

\item[] \textsc{Goldfarb, D.} and \textsc{Iyengar, G.} (2003). Robust portfolio selection problems.
{\em Math. Oper. Res.} {\bf 28} 1--38.

\item[] \textsc{Horn, R. A.} and \textsc{Johnson, C. R.} (1990). {\em Matrix Analysis}.
Cambridge University Press, Cambridge.

\item[] \textsc{Huang, J. Z.}, \textsc{Liu, N.} and \textsc{Pourahmadi, M.} (2004). Covariance selection
and esimation via penalized normal likelihood. {\em Manuscript}.

\item[] \textsc{Huber, P.} (1973). Robust regression: asymptotics, conjectures and
Monte Carlo. {\em Ann. Statist.} {\bf 1} 799--821.

\item[] \textsc{James, W.} and \textsc{Stein, C.} (1961). Estimation with quadratic
Loss. In {\em Proc. Fourth Berkeley Symp. Math. Statist. Probab.}
{\bf 1} 361--379. Univ. California Press, Berkeley.

\item[] \textsc{Johnstone, I. M.} (2001). On the distribution of the largest eigenvalue in
principal components analysis. {\em Ann. Statist.} {\bf 29}
295--327.

\item[] \textsc{Ledoit, O.} and \textsc{Wolf, M.} (2004). A well conditioned estimator for
large--dimensional covariance matrices. {\em J. Multiv. Anal.} {\bf
88} 365--411.

\item[] \textsc{Leonard, T.} and \textsc{Hsu, J.S.J.} (1992). Bayesian inference for a covariance
matrix. {\em Ann. Statist.} {\bf 20} 1669--1696.

\item[] \textsc{Li, H.} and \textsc{Gui, J.} (2005). Gradient directed regularization for sparse Gaussian
concentration graphs, with applications to inference of genetic
networks. {\em Manuscript}.

\item[] \textsc{Lin, S. P.} and \textsc{Perlman, M. D.} (1985). A Monte Carlo comparison of four estimators of
a covariance matrix. In {\em Multivariate Analysis} {\bf 6}, Ed. P.
R. Krishnaiah, 411--429. Amsterdam, North-Holland.

\item[] \textsc{Markowitz, H. M.} (1952). Portfolio selection. {\em Journal of
Finance} {\bf 7} 77--91.

\item[] \textsc{Markowitz, H. M.} (1959). {\em Portfolio Selection: Efficient Diversification
of Investments}. John Wiley \& Sons, New Jersey.

\item[] \textsc{Muirhead, Robb J.} (1982). {\em Aspects of Multivariate Statistical
Theory}. John Wiley \& Sons, New York.

\item[] \textsc{Portnoy, S.} (1984). Asymptotic behavior of M-estimators of $p$ regression
parameters when $p^2/n$ is large. I. Consistency. {\em Ann.
Statist.} {\bf 12} 1298--1309.

\item[] \textsc{Portnoy, S.} (1985). Asymptotic behavior of M estimators of $p$ regression
parameters when $p^2/n$ is large. II. Normal approximation. {\em
Ann. Statist.} {\bf 13} 1403--1417.

\item[] \textsc{Pourahmadi, M.} (2000). Maximum likelihood estimation of generalized linear
models for multivariate normal covariance matrix. {\em Biometrika}
{\bf 87} 425--435.

\item[] \textsc{Ross, S. A.} (1976). The Arbitrage Theory of Capital Asset
Pricing. {\em Journal of Economic Theory} {\bf 13} 341--360.

\item[] -------- (1977). The Capital Asset Pricing Model (CAPM), short-sale restrictions
and related issues. {\em Journal of Finance} {\bf 32} 177--183.

\item[] \textsc{Scott, J. T.} (1966). Factor analysis and regression. {\em Econometrica} {\bf 34}
552--562.

\item[] \textsc{Scott, J. T.} (1969). Factor analysis regression revisited. {\em Econometrica}
{\bf 37} 719.

\item[] \textsc{Smith, M.} and \textsc{Kohn, R.} (2002). Parsimonious covariance matrix
estimation for longitudinal data. {\em J. Amer. Statist. Assoc.}
{\bf 97} 1141--1153.

\item[] \textsc{Stein, C.} (1975). Estimation of a covariance matrix. {\em Rietz
Lecture}, 39th IMS Annual Meeting, Atlanta, Georgia.

\item[] \textsc{Stock, J. H.} and \textsc{Watson, M. W.} (2005). Implications of dynamic factor models for VAR
analysis. {\em Manuscript}.

\item[] \textsc{Wong, F.}, \textsc{Carter, C. K.} and \textsc{Kohn, R.} (2003). Efficient estimation of
covariance selection models. {\em Biometrika} {\bf 90} 809--830.

\item[] \textsc{Wu, W. B.} and \textsc{Pourahmadi, M.} (2003). Nonparametric estimation of large covariance
matrices of longitudinal data. {\em Biometrika} {\bf 90} 831--844.

\item[] \textsc{Yohai, V. J.} and \textsc{Maronna, R. A.} (1979). Asymptotic behavior of M-estimators
for the linear model. {\em Ann. Statist.} {\bf 7} 258--268.

\item[] \textsc{Yuan, K.-H.} and \textsc{Bentler, P. M.} (1997). Mean and covariance
structure analysis: theoretical and practical improvements. {\em J.
Amer. Statist. Assoc.} {\bf 438} 767--774.
\end{itemize}
\end{singlespace}

\bigskip


\footnotesize
\begin{quote}
\begin{tabular}{lcl}
\textsc{Jianqing Fan} & \quad\quad\quad & \textsc{Yingying Fan}\\
\textsc{Department of Operations Research} & \quad\quad\quad & \textsc{Department of Operations Research}\\
\textsc{\ \ and Financial Engineering} & & \textsc{\ \ and Financial Engineering}\\
\textsc{Princeton University} & & \textsc{Princeton
University}\\
\textsc{Princeton, New Jersey 08544} & & Princeton, New Jersey 08544\\
USA & & USA\\
\textsc{E-mail}: jqfan@princeton.edu & & \textsc{E-mail}:
yingying@princeton.edu
\end{tabular}
\\
\begin{center}
\begin{tabular}{l}
\textsc{Jinchi Lv}\\
\textsc{Department of Mathematics}\\
\textsc{Princeton University}\\
\textsc{Princeton, New Jersey
08544}\\
USA\\
\textsc{E-mail}: jlv@princeton.edu\\
\end{tabular}
\end{center}
\end{quote}

\end{document}